\documentclass{article}

\usepackage{amsfonts,amsmath,amssymb}
\usepackage{graphicx}
\usepackage{color}
\usepackage{subcaption}
\usepackage{placeins}
\usepackage{hyperref}
\usepackage{tikz}

\definecolor{darkblue}{rgb}{0,0.08,0.45}
\hypersetup{colorlinks,linkcolor=darkblue,citecolor=darkblue,filecolor=darkblue,urlcolor=darkblue}
\hypersetup{bookmarksopen,bookmarksnumbered}

\newcommand{\ipg}[2]{\ensuremath{\left\langle #1 , #2 \right\rangle_\Gamma}}
\newcommand{\norm}[1]{\ensuremath{\left\lVert #1 \right\rVert}}
\newcommand{\normg}[1]{\ensuremath{\left\lVert #1 \right\rVert_\Gamma}}
\newcommand{\dx}{\ensuremath{\,\mathrm{d}\mathbf{x}}}
\newcommand{\dy}{\ensuremath{\,\mathrm{d}\mathbf{y}}}
\newcommand{\normal}{\ensuremath{\hat{\mathbf{n}}}}
\newcommand{\ptot}{\ensuremath{p_{\mathrm{tot}}}}
\newcommand{\psca}{\ensuremath{p_{\mathrm{sca}}}}
\newcommand{\pinc}{\ensuremath{p_{\mathrm{inc}}}}
\newcommand{\pint}{\ensuremath{p_{\mathrm{int}}}}
\newcommand{\SLP}{\ensuremath{\mathcal{V}}}
\newcommand{\DLP}{\ensuremath{\mathcal{K}}}
\newcommand{\SL}{\ensuremath{V}}
\newcommand{\DL}{\ensuremath{K}}
\newcommand{\AD}{\ensuremath{T}}
\newcommand{\HS}{\ensuremath{D}}
\newcommand{\ID}{\ensuremath{I}}
\newcommand{\REG}{\ensuremath{R}}
\newcommand{\REGntd}{\ensuremath{\REG_\mathrm{NtD}}}
\newcommand{\TR}{\ensuremath{Z}}
\newcommand{\TRc}{\ensuremath{\overline{\TR}}}
\newcommand{\LB}{\ensuremath{\Delta_\Gamma}}
\newcommand{\NtDe}{\ensuremath{\Lambda_\mathrm{NtD}^+}}
\newcommand{\DtNe}{\ensuremath{\Lambda_\mathrm{DtN}^+}}
\newcommand{\OsrcNtD}{\ensuremath{L_\mathrm{NtD}}}
\newcommand{\OsrcDtN}{\ensuremath{L_\mathrm{DtN}}}
\newcommand{\Hminushalf}{\ensuremath{H^{-1/2}(\Gamma)}}
\newcommand{\Hplushalf}{\ensuremath{H^{1/2}(\Gamma)}}
\newcommand{\kint}{\ensuremath{k_\text{int}}}
\newcommand{\kext}{\ensuremath{k_\text{ext}}}
\newcommand{\rhoint}{\ensuremath{\rho_\text{int}}}
\newcommand{\rhoext}{\ensuremath{\rho_\text{ext}}}
\newcommand{\cint}{\ensuremath{c_\text{int}}}
\newcommand{\cext}{\ensuremath{c_\text{ext}}}
\newcommand{\traceDe}{\ensuremath{\gamma_D^+}}
\newcommand{\traceDi}{\ensuremath{\gamma_D^-}}
\newcommand{\traceDei}{\ensuremath{\gamma_D^\pm}}
\newcommand{\traceNe}{\ensuremath{\gamma_N^+}}
\newcommand{\traceNi}{\ensuremath{\gamma_N^-}}
\newcommand{\traceNei}{\ensuremath{\gamma_N^\pm}}

\title{Stable and efficient FEM-BEM coupling with OSRC regularisation for acoustic wave transmission\footnote{© 2021. This manuscript version is made available under the CC-BY-NC-ND 4.0 license. This manuscript is published in the Journal of Computational Physics in final form at \url{https://doi.org/10.1016/j.jcp.2021.110867}.}}
\author{Elwin van 't Wout\thanks{Institute for Mathematical and Computational Engineering, School of Engineering and Faculty of Mathematics, Pontificia Universidad Católica de Chile, Santiago, Chile. Contact: e.wout@uc.cl}}
\date{November 26, 2021}

\begin{document}
	
\maketitle

\begin{abstract}
	The finite element method (FEM) and the boundary element method (BEM) can numerically solve the Helmholtz system for acoustic wave propagation. When an object with heterogeneous wave speed or density is embedded in an unbounded exterior medium, the coupled FEM-BEM algorithm promises to combine the strengths of each technique. The FEM handles the heterogeneous regions while the BEM models the homogeneous exterior. Even though standard FEM-BEM algorithms are effective, they do require stabilisation at resonance frequencies. One such approach is to add a regularisation term to the system of equations. This algorithm is stable at all frequencies but also brings higher computational costs. This study proposes a regulariser based on the on-surface radiation conditions (OSRC). The OSRC operators are also used to precondition the boundary integral operators and combined with incomplete LU factorisations for the volumetric weak formulation. The proposed preconditioning strategy improves the convergence of iterative linear solvers significantly, especially at higher frequencies.
\end{abstract}

\section{Introduction}

The Helmholtz equation models a wide range of propagation phenomena for harmonic acoustic waves~\cite{lahaye2017modern}. In acoustic transmission problems, the model involves a system of partial differential equations that can be solved numerically with the finite element method (FEM) or the boundary element method (BEM), among other discretisation techniques. When the model includes acoustic scattering in an unbounded exterior domain, the surface-based BEM is the preferred option since the Sommerfeld radiation condition is automatically satisfied by the representation formula. However, when heterogeneous structures are present, the absence of a Green's function prohibits the BEM and volumetric methods like the FEM are required. Hence, many techniques have been proposed in the literature to couple the FEM and BEM, combining the complementary strengths of both methodologies and simulate acoustic transmission into heterogeneous materials embedded in an unbounded medium (cf.~\cite{johnson1980coupling, costabel1987symmetric, hiptmair2006stabilized, stephan2018coupling}).

This study considers strong coupling of the FEM and BEM, where a single model explicitly includes the wave propagation in the entire geometry. The system matrix includes blocks corresponding to the FEM for the bounded heterogeneous subdomains, the BEM for the unbounded exterior medium, and coupling operators at the material interface. An alternative approach is weak coupling, where the FEM and BEM models are solved independently in their respective domains and then iteratively coupled with the transmission conditions. On the one hand, weakly coupled FEM-BEM algorithms are easier to implement since separate software for the FEM and BEM can be used~\cite{caudron2020optimized}. On the other hand, strongly coupled FEM-BEM algorithms are more accurate~\cite{gaul2008coupling}. Other approaches of coupling the BEM with volumetric methods include mortar element methods~\cite{mascotto2020fem}, tearing-and-interconnection techniques~\cite{langer2005coupled}, domain decomposition~\cite{langer2007coupled}, and surface-volume integral equations~\cite{costabel2015spectrum}. The coupling techniques are available for the time-domain wave equation as well~\cite{abboud2011coupling}. Finally, one can also interpret the coupled FEM-BEM algorithm as a model where the BEM is applied as an artificial boundary condition to truncate the computational domain of the FEM.

A drawback of most FEM-BEM algorithms is their instability at resonance frequencies of the geometry. Specialised stabilisation techniques need to be employed to avoid spurious solutions and ill-conditioned systems. One such stabilisation approach introduces a regularisation variable in the FEM-BEM system~\cite{hiptmair2006stabilized}. This stabilised FEM-BEM algorithm can be proven to be stable at all frequencies when the regularisation operator is positive-definite. Several examples of such regularisers were proposed in literature, including local boundary integral operators~\cite{buffa2005regularized}. Furthermore, the effectiveness of this stabilisation approach has been shown for the convected Helmholtz equations~\cite{casenave2014coupled}, elastodynamics~\cite{gatica2014coupling}, and electromagnetics~\cite{hiptmair2008stabilized} as well. However, this regularisation technique incurs higher computational costs, also at non-resonant frequencies. Specifically, the stabilised system includes more degrees of freedom and is worse conditioned compared to standard formulations.

This study deals with the computational efficiency of the stabilised FEM-BEM algorithm. It introduces a novel regularisation operator for the stabilisation term and preconditioning of the system matrix: an on-surface radiation condition (OSRC) approximation of the Neumann-to-Dirichlet (NtD) map. The OSRC operators have been used in a variety of acoustic algorithms, including scattering models~\cite{kriegsmann1987new}, radiation conditions~\cite{tsynkov1998numerical}, coupling operators~\cite{antoine2004analytic}, and preconditioners~\cite{wout2021pmchwt}. They are especially effective at high frequencies~\cite{antoine2021introduction}. In the context of FEM-BEM coupled systems, the OSRC operator has been used for weak coupling~\cite{caudron2020optimized}, but not yet for strong coupling.

Solving the system matrix with iterative linear algorithms can drastically be improved with preconditioning. The design of efficient preconditioners is highly problem-dependent and is challenging for FEM-BEM coupling since the FEM and BEM operators have different properties~\cite{stephan2018coupling, feischl2017optimal}. This study uses operator preconditioning with OSRC operators~\cite{antoine2021introduction}. Various large-scale simulations of high-frequency wave propagation have already confirmed the efficiency of OSRC preconditioning for pure BEM~\cite{haqshenas2021fast, wout2021pmchwt}. Here, this strategy is adapted to the FEM-BEM algorithm and combined with incomplete LU factorisation for the volumetric part of the weak formulation. The computational benchmarks in this study show that preconditioning and stabilisation with OSRC operators significantly improve the efficiency of the stabilised FEM-BEM algorithm.

Section~\ref{sec:formulation} presents the formulation of the OSRC-stabilised FEM-BEM system and the preconditioners. The computational results in Section~\ref{sec:results} confirm the stability and efficiency of the proposed methodology.

\section{Formulation}
\label{sec:formulation}

This section summarises the Helmholtz system for acoustic transmission and presents the stable coupling of the FEM and BEM with OSRC regularisation, followed by the design of efficient preconditioners.

\subsection{Model equations}

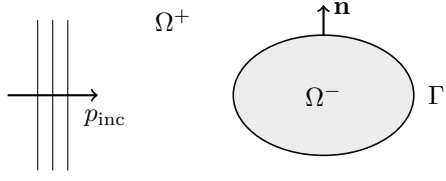
\begin{figure}[!ht]
	\centering
	\begin{tikzpicture}
		\draw[thin](-4.8,-1)--(-4.8,1);
		\draw[thin](-4.6,-1)--(-4.6,1);
		\draw[thin](-4.4,-1)--(-4.4,1);
		\draw[thick, ->](-5.2,0)--(-4,0);
		\node at (-3.9,-0.3) {$\pinc$};
		
		\node at (-3,1) {$\Omega^+$};
		
		\draw[semithick, fill=gray!15](-1,0) ellipse (1.2 and 0.8) node {$\Omega^-$};
		\node at (0.5,0) {$\Gamma$};
		\draw[thick, ->](-1,0.8)--(-1,1.2) node[right] {$\normal$};
	\end{tikzpicture}
	\caption{A sketch of the geometry.}
	\label{fig:domain}
\end{figure}

Let us consider a structure embedded in free space, as depicted in Fig.~\ref{fig:domain}. The volume of the bounded object is denoted by $\Omega^- \subset \mathbb{R}^3$ and has a piecewise smooth boundary $\Gamma$ which can be equipped with a unit outward-pointing normal $\normal$. The exterior domain, denoted by $\Omega^+$, is the complement and assumed to be homogeneous with constant density $\rhoext$ and wave speed $\cext$. The interior domain can be heterogeneous with density $\rhoint(\mathbf{x})$ and wave speed $\cint(\mathbf{x})$ for $\mathbf{x} \in \Omega^-$. No attenuation is considered so that the wavenumber satisfies
\begin{equation}
	k(\mathbf{x}) = \frac{2\pi f}{c(\mathbf{x})}
\end{equation}
where $f$ denotes the frequency of the known source field $\pinc$. The density and wave speed are assumed to vary smoothly inside the volume $\Omega^-$ but can have a jump across the surface $\Gamma$.

The acoustic pressure field in the exterior is denoted by $\ptot$ and the scattered field by $\psca = \ptot - \pinc$. The acoustic pressure field in the interior is denoted by $\pint$. Harmonic wave propagation through a material with a linear response can be modelled by the Helmholtz equation as
\begin{equation}
	\label{eq:helmholtz}
	\begin{cases}
		-\nabla^2 \psca - \kext^2 \psca = 0 & \text{in } \Omega^+, \\
		-\rhoint \nabla \cdot \left( \frac1\rhoint \nabla \pint \right) - \kint^2 \pint = 0 & \text{in } \Omega^-, \\
		\traceDe \ptot = \traceDi \pint & \text{on } \Gamma, \\
		\traceNe \left(\frac1\rhoext \ptot\right) = \traceNi \left(\frac1\rhoint \pint\right) & \text{on } \Gamma, \\
		\lim_{\mathbf{x} \to \infty} |\mathbf{x}| \left(\partial_{|\mathbf{x}|} \psca - \imath\kext \psca\right) = 0
	\end{cases}
\end{equation}
where the last equation is the Sommerfeld radiation condition, the Dirichlet and Neumann traces are defined by
\begin{subequations}
\begin{align}
	\traceDei p(\mathbf{x}) &= \lim_{\mathbf{y} \to \mathbf{x}} p(\mathbf{y}) && \text{for } \mathbf{x} \in \Gamma \text{ and } \mathbf{y} \in \Omega^\pm, \\
	\traceNei p(\mathbf{x}) &= \lim_{\mathbf{y} \to \mathbf{x}} \nabla p(\mathbf{y}) \cdot \normal && \text{for } \mathbf{x} \in \Gamma \text{ and } \mathbf{y} \in \Omega^\pm
\end{align}
\end{subequations}
where the $\pm$ denote the exterior and interior versions, respectively, and $\imath$ denotes the complex unit.

\subsection{Exterior boundary integral formulation}

The pressure field in the unbounded exterior domain will be represented by surface potentials. Specifically, any radiating solution of the Helmholtz system can be written as
\begin{equation}
	\label{eq:representation}
	\psca = \DLP(\traceDe\ptot) - \SLP(\traceNe\ptot)
\end{equation}
with $\SLP$ and $\DLP$ defined by
\begin{subequations}
	\label{eq:potentialoperators}
	\begin{align}
		[\SLP\psi](\mathbf{x}) &= \iint_\Gamma G(\mathbf{x},\mathbf{y}) \psi(\mathbf{y}) \dy && \text{for } \mathbf{x} \in \Omega^+, \\
		[\DLP\phi](\mathbf{x}) &= \iint_\Gamma \frac{\partial G(\mathbf{x},\mathbf{y})}{\partial \normal(\mathbf{y})} \phi(\mathbf{y}) \dy && \text{for } \mathbf{x} \in \Omega^+,
	\end{align}
\end{subequations}
the single-layer and double-layer potential operators, respectively, and
\begin{equation}
	\label{eq:green}
	G(\mathbf{x},\mathbf{y}) = \frac{e^{\imath\kext|\mathbf{x}-\mathbf{y}|}}{4\pi|\mathbf{x}-\mathbf{y}|}, \qquad \mathbf{x},\mathbf{y} \in \Omega^+, \mathbf{x} \ne \mathbf{y}
\end{equation}
the Green's function (cf.~\cite{nedelec2001acoustic, steinbach2008numerical, sauter2010boundary}).
As shown in~\ref{sec:calderon}, given a valid integral operator
\begin{equation*}
	\mathbb{P} = \begin{bmatrix} \mathbb{A} & \mathbb{B} \\ \mathbb{C} & \mathbb{D} \end{bmatrix}
\end{equation*}
that models the Helmholtz system (c.f.~Eq.~\eqref{eq:calderon:projector}), the traces of the scattered field satisfy the Calderón projection property
\begin{equation}
	\begin{bmatrix} \traceDe\psca \\ \traceNe\psca \end{bmatrix} = \mathbb{P} \begin{bmatrix} \traceDe\psca \\ \traceNe\psca \end{bmatrix}.
\end{equation}
The Dirichlet-to-Neumann (DtN) and Neumann-to-Dirichlet (NtD) maps are implicitly defined by
\begin{subequations}
\begin{align}
	\DtNe\traceDe\psca &= \traceNe\psca, \\
	\NtDe\traceNe\psca &= \traceDe\psca
\end{align}
\end{subequations}
so that the DtN map is the inverse of the NtD map. Different expressions of the DtN map exist, including
\begin{subequations}
\begin{align}
	\DtNe &= \mathbb{B}^{-1} (\ID - \mathbb{A}), \\
	\DtNe &= \mathbb{C} + \mathbb{D} \mathbb{B}^{-1} \left(\ID - \mathbb{A}\right),
\end{align}
\end{subequations}
as derived in~\ref{sec:dtn}. Since the DtN and NtD maps include inverse operators, no explicit expressions are available for general surfaces.

The Calderón projector for the standard boundary integral formulations reads
\begin{align}
	\mathbb{P} &= \begin{bmatrix} \tfrac12\ID + \DL & -\SL \\ -\HS & \tfrac12\ID - \AD \end{bmatrix}
\end{align}
and the DtN map satisfies
\begin{subequations}
\begin{align}
	\DtNe &= \SL^{-1} \left(\DL - \tfrac12\ID\right), \label{eq:dtn:standard} \\
	\DtNe &= -\HS + \left(\tfrac12\ID - \AD\right) \SL^{-1} \left(\DL - \tfrac12\ID\right).
\end{align}
\end{subequations}
Here, $\ID$ denotes the identity operator and
\begin{subequations}
\label{eq:boundaryoperators}
\begin{align}
	[\SL\psi](\mathbf{x}) &= \iint_\Gamma G(\mathbf{x},\mathbf{y}) \psi(\mathbf{y}) \dy && \text{for } \mathbf{x} \in \Gamma, \\
	[\DL\phi](\mathbf{x}) &= \iint_\Gamma \frac{\partial}{\partial \normal(\mathbf{y})} G(\mathbf{x},\mathbf{y}) \phi(\mathbf{y}) \dy && \text{for } \mathbf{x} \in \Gamma, \\
	[\AD\psi](\mathbf{x}) &= \frac{\partial}{\partial \normal(\mathbf{x})} \iint_\Gamma G(\mathbf{x},\mathbf{y}) \psi(\mathbf{y}) \dy && \text{for } \mathbf{x} \in \Gamma, \\
	[\HS\phi](\mathbf{x}) &= -\frac{\partial}{\partial \normal(\mathbf{x})} \iint_\Gamma \frac{\partial}{\partial \normal(\mathbf{y})} G(\mathbf{x},\mathbf{y}) \phi(\mathbf{y}) \dy && \text{for } \mathbf{x} \in \Gamma
\end{align}
\end{subequations}
are called the single-layer, double-layer, adjoint double-layer and hypersingular boundary integral operator, respectively.

\subsection{Interior weak formulation}

In general, wave propagation through a heterogeneous material does not satisfy a Green's function and the BEM cannot be used anymore. Hence, the FEM will be used for the interior domain. The weak formulation of the Helmholtz equation is to search for a solution $\pint \in H^1(\Omega^-)$ such that
\begin{equation}
	\mathsf{a}(\pint,q) - \ipg{\traceNi\pint}{\traceDi q} = 0, \quad \forall q \in H^1(\Omega^-)
\end{equation}
with the sesquilinear form
\begin{align}
	\mathsf{a}(\pint,q) &= \iiint_{\Omega^-} \frac1{\rhoint(\mathbf{x})} \nabla \pint(\mathbf{x}) \nabla \left(\rhoint(\mathbf{x}) q(\mathbf{x})\right) \dx \nonumber \\ &\qquad - \iiint_{\Omega^-} \kint^2(\mathbf{x}) \pint(\mathbf{x}) q(\mathbf{x}) \dx
\end{align}
where $\ipg{\cdot}{\cdot}$ denotes the standard $L_2$ inner product on $\Gamma$ and $H^1$ the standard Sobolev space on $\Omega^-$ (cf.~\cite{steinbach2008numerical}). To obtain a consistent set of equations, the Neumann trace of the interior field needs to be written in terms of the Dirichlet trace. This can be achieved with the DtN map as follows:
\begin{align*}
	&\ipg{\traceNi\pint}{\traceDi q} \\
	&= \ipg{\tfrac\rhoext\rhoint \traceNe\ptot}{\traceDi q} \\
	&= \ipg{\traceNe\psca + \traceNe\pinc}{\tfrac\rhoext\rhoint \traceDi q} \\
	&= \ipg{\DtNe \traceDe\psca}{\tfrac\rhoext\rhoint \traceDi q} + \ipg{\traceNe\pinc}{\tfrac\rhoext\rhoint \traceDi q} \\
	&= \ipg{\left(\mathbb{C} + \mathbb{D} \mathbb{B}^{-1} \left(\ID - \mathbb{A}\right)\right) \traceDe(\ptot - \pinc)}{\tfrac\rhoext\rhoint \traceDi q} + \ipg{\traceNe\pinc}{\tfrac\rhoext\rhoint \traceDi q} \\
	&= \ipg{\mathbb{C} (\traceDe\ptot - \traceDe\pinc)}{\tfrac\rhoext\rhoint \traceDi q} + \ipg{\mathbb{D} \vartheta}{\tfrac\rhoext\rhoint \traceDi q} + \ipg{\traceNe\pinc}{\tfrac\rhoext\rhoint \traceDi q} \\
	&= \ipg{\mathbb{C} \traceDi\pint}{\tfrac\rhoext\rhoint \traceDi q} + \ipg{\mathbb{D} \vartheta}{\tfrac\rhoext\rhoint \traceDi q} - \ipg{\mathbb{C} \traceDe\pinc}{\tfrac\rhoext\rhoint \traceDi q} + \ipg{\traceNe\pinc}{\tfrac\rhoext\rhoint \traceDi q}
\end{align*}
for $\vartheta$ an unknown surface potential that satisfies
\begin{equation}
	\mathbb{B} \vartheta = \left(\ID - \mathbb{A}\right) \traceDi\pint + \left(\mathbb{A} - \ID\right) \traceDe\pinc.
\end{equation}
The coupled variational problem is now given by: search for a solution $\pint \in H^1(\Omega^-)$ and surface potential $\vartheta \in H^{-1/2}(\Gamma)$ such that
\begin{subequations}
\label{eq:fembem:symmetric:general}
\begin{align}
	&\mathsf{a}(\pint,q) - \ipg{\tfrac\rhoext\rhoint \mathbb{C} \traceDi\pint}{\traceDi q} - \ipg{\tfrac\rhoext\rhoint \mathbb{D} \vartheta}{\traceDi q} \nonumber \\ &\qquad = -\ipg{\tfrac\rhoext\rhoint \mathbb{C} \traceDe\pinc}{\traceDi q} + \ipg{\tfrac\rhoext\rhoint \traceNe\pinc}{\traceDi q}, && \forall q \in H^1(\Omega^-); \\
	&\ipg{\varphi}{\left(\ID - \mathbb{A}\right) \traceDi\pint} - \ipg{\varphi}{\mathbb{B} \vartheta} \nonumber \\ &\qquad = \ipg{\varphi}{\left(\ID - \mathbb{A}\right) \traceDe\pinc}, && \forall \varphi \in H^{-1/2}(\Gamma).
\end{align}
\end{subequations}
Substituting the standard Calderón operators yields: search for a solution $\pint \in H^1(\Omega^-)$ and surface potential $\vartheta \in H^{-1/2}(\Gamma)$ such that
\begin{subequations}
\label{eq:fembem:symmetric}
\begin{align}
	&\mathsf{a}(\pint,q) + \ipg{\tfrac\rhoext\rhoint \HS \traceDi\pint}{\traceDi q} + \ipg{\tfrac\rhoext\rhoint (\AD - \tfrac12\ID) \vartheta}{\traceDi q} \nonumber \\ &\qquad = \ipg{\tfrac\rhoext\rhoint \HS \traceDe\pinc}{\traceDi q} + \ipg{\tfrac\rhoext\rhoint \traceNe\pinc}{\traceDi q}, && \forall q \in H^1(\Omega^-); \\
	&\ipg{\varphi}{\left(\tfrac12\ID - \DL\right) \traceDi\pint} + \ipg{\varphi}{\SL \vartheta} \nonumber \\ &\qquad = \ipg{\varphi}{\left(\tfrac12\ID - \DL\right) \traceDe\pinc}, && \forall \varphi \in H^{-1/2}(\Gamma).
\end{align}
\end{subequations}
This variational problem is called the symmetric FEM-BEM coupling algorithm~\cite{costabel1987symmetric} and is known to be unstable at the resonance frequencies of the geometry~\cite{hiptmair2006stabilized}.

\subsection{Stabilised FEM-BEM coupling}

The idea behind the stabilised FEM-BEM coupling in~\cite{hiptmair2006stabilized} is to transform the solution traces such that stability at resonance frequencies can be proven. Specifically, defining a general trace transformation operator as
\begin{equation}
	\mathcal{T} = \begin{bmatrix} \mathcal{A} & \mathcal{B} \\ \mathcal{C} & \mathcal{D} \end{bmatrix}
\end{equation}
one can use the Calderón projector
\begin{equation}
	\mathbb{P} = \begin{bmatrix} \mathcal{A}^{-1} & 0 \\ 0 & \mathcal{D}^{-1} \end{bmatrix} \left( \begin{bmatrix} \mathcal{A} & \mathcal{B} \\ \mathcal{C} & \mathcal{D} \end{bmatrix} \begin{bmatrix} \tfrac12\ID + \DL & -\SL \\ -\HS & \tfrac12\ID - \AD \end{bmatrix} - \begin{bmatrix} 0 & \mathcal{B} \\ \mathcal{C} & 0 \end{bmatrix} \right)
\end{equation}
that will solve the Helmholtz transmission problem (c.f.~\cite[Lemma~6.3]{hiptmair2006stabilized}).
A convenient choice of valid trace transformation is given by
\begin{equation}
	\label{eq:tracetransformation}
	\mathcal{T}_\text{reg} = \begin{bmatrix} \ID & \imath\eta\REG \\ \imath\nu\ID & \ID \end{bmatrix}
\end{equation}
where $\eta \in \mathbb{R}$, $\eta \ne 0$, $\nu=0$ or $\nu=\eta$, and
\begin{equation}
	\label{eq:def:reg}
	\REG: H^{-1/2}(\Gamma) \to H^{1/2}(\Gamma)
\end{equation}
a regularising operator (c.f.~\cite[Eq.~(39)]{hiptmair2006stabilized}).
The choice of the regulariser $\REG$ needs to satisfy
\begin{enumerate}
	\item $\REG: H^{-1/2}(\Gamma) \to H^{1/2}(\Gamma)$ is compact, and
	\item $\Re\left(\ipg{\theta}{\REG(\theta)}\right) > 0$ for all $\theta \in H^{-1/2}(\Gamma)\setminus\{0\}$
\end{enumerate}
to prove stability of the FEM-BEM formulation at all frequencies~\cite{hiptmair2006stabilized, meury2007stable, buffa2005regularized, gatica2014coupling}.
Notice that the standard Calderón projector is retrieved for $\eta=0$.
Now,
\begin{subequations}
\begin{align}
	\mathbb{A} &= \tfrac12\ID + \DL - \mathcal{A}^{-1} \mathcal{B} \HS = \tfrac12\ID + \DL - \imath\eta\REG \HS, \\
	\mathbb{B} &= -\SL - \mathcal{A}^{-1} \mathcal{B} (\tfrac12\ID + \AD) = -\SL - \imath\eta\REG (\tfrac12\ID + \AD), \\
	\mathbb{C} &= -\mathcal{D}^{-1} \mathcal{C} (\tfrac12\ID - \DL) - \HS = -\imath\nu (\tfrac12\ID - \DL) - \HS, \\
	\mathbb{D} &= -\mathcal{D}^{-1} \mathcal{C} \SL + \tfrac12\ID - \AD = -\imath\nu \SL + \tfrac12\ID - \AD.
\end{align}
\end{subequations}
Substitution into the symmetric FEM-BEM formulation~\eqref{eq:fembem:symmetric:general} yields
\begin{align*}
	&\mathsf{a}(\pint,q) + \ipg{\tfrac\rhoext\rhoint \left(\imath\nu (\tfrac12\ID - \DL) + \HS\right) \traceDi\pint}{\traceDi q} + \ipg{\tfrac\rhoext\rhoint \left(\imath\nu \SL + \AD - \tfrac12\ID\right) \vartheta}{\traceDi q} \nonumber \\ &\qquad = \ipg{\tfrac\rhoext\rhoint \left(\imath\nu (\tfrac12\ID - \DL) + \HS\right) \traceDe\pinc}{\traceDi q} + \ipg{\tfrac\rhoext\rhoint \traceNe\pinc}{\traceDi q}, \\
	&\ipg{\varphi}{\left(\tfrac12\ID - \DL + \imath\eta\REG \HS\right) \traceDi\pint} + \ipg{\varphi}{\left(\SL + \imath\eta\REG (\tfrac12\ID + \AD)\right) \vartheta} \nonumber \\ &\qquad = \ipg{\varphi}{\left(\tfrac12\ID - \DL + \imath\eta\REG \HS\right) \traceDe\pinc}.
\end{align*}
As will be shown in Section~\ref{sec:regularisation}, many regularising operators are defined as the inverse of a differential operator. Hence, let us introduce a new variable
\begin{equation}
	\sigma = \REG (\tfrac12\ID + \AD) \vartheta + \REG \HS (\traceDi\pint - \traceDe\pinc)
\end{equation}
which reduces the second equation in the FEM-BEM system to
\begin{align*}
	&\ipg{\varphi}{\left(\tfrac12\ID - \DL\right) \traceDi\pint} + \ipg{\varphi}{\imath\eta\REG \HS \traceDi\pint} + \ipg{\varphi}{\SL \vartheta} + \ipg{\varphi}{\imath\eta\REG (\tfrac12\ID + \AD) \vartheta} \nonumber \\ &\quad = \ipg{\varphi}{\left(\tfrac12\ID - \DL\right) \traceDe\pinc} + \ipg{\varphi}{\imath\eta\REG \HS \traceDe\pinc}, \\
	&\ipg{\varphi}{\left(\tfrac12\ID - \DL\right) \traceDi\pint} + \ipg{\varphi}{\SL \vartheta} + \ipg{\varphi}{\imath\eta\ID \sigma} = \ipg{\varphi}{\left(\tfrac12\ID - \DL\right) \traceDe\pinc}.
\end{align*}
The newly introduced variable $\sigma$ needs to satisfy the variational formulation
\begin{equation}
	\mathsf{b}(\sigma,\tau) - \ipg{(\tfrac12\ID + \AD) \vartheta}{\tau} - \ipg{\HS \traceDi\pint}{\tau} = -\ipg{\HS \traceDe\pinc}{\tau}
\end{equation}
where
\begin{equation}
	\label{eq:b}
	\mathsf{b}(\REG(\sigma),\tau) = \ipg{\sigma}{\tau}
\end{equation}
a sesquilinear form that represents the weak form of the inverse regulariser: $\mathsf{b}(\sigma,\tau) = \ipg{\REG^{-1}\sigma}{\tau}$. Putting everything together, the variational problem reads: search for $\pint \in H^1(\Omega^-)$, $\vartheta \in H^{-1/2}(\Gamma)$ and $\sigma \in H^1(\Gamma)$ such that
\begin{subequations}
\label{eq:fembem:stable}
\begin{align}
	&\mathsf{a}(\pint,q) + \ipg{\tfrac\rhoext\rhoint \left(\imath\nu (\tfrac12\ID - \DL) + \HS\right) \traceDi\pint}{\traceDi q} \nonumber \\ &\qquad + \ipg{\tfrac\rhoext\rhoint \left(\imath\nu \SL + \AD - \tfrac12\ID\right) \vartheta}{\traceDi q} \nonumber \\ &\qquad = \ipg{\tfrac\rhoext\rhoint \left(\imath\nu (\tfrac12\ID - \DL) + \HS\right) \traceDe\pinc}{\traceDi q} \nonumber \\ &\qquad\qquad + \ipg{\tfrac\rhoext\rhoint \traceNe\pinc}{\traceDi q}, && \forall q \in H^1(\Omega^-); \\
	&\ipg{\varphi}{\left(\tfrac12\ID - \DL\right) \traceDi\pint} + \ipg{\varphi}{\SL \vartheta} + \ipg{\varphi}{\imath\eta\ID \sigma} \nonumber \\ &\qquad = \ipg{\varphi}{\left(\tfrac12\ID - \DL\right) \traceDe\pinc}, && \forall \varphi \in H^{-1/2}(\Gamma); \\
	&\mathsf{b}(\sigma,\tau) - \ipg{(\tfrac12\ID + \AD) \vartheta}{\tau} - \ipg{\HS \traceDi\pint}{\tau} \nonumber \\ &\qquad = -\ipg{\HS \traceDe\pinc}{\tau}, && \forall \tau \in H^1(\Gamma),
\end{align}
\end{subequations}
which is the stabilised FEM-BEM formulation~\cite[Eq.~(46)]{hiptmair2006stabilized}.

\subsubsection{Interpretation of variables}

The variational formulation of the stabilised FEM-BEM algorithm includes three unknowns. Firstly, $\pint$ is the pressure field inside the object. Secondly, $\vartheta$ is a surface potential that satisfies
\begin{align*}
	\vartheta
	&= \mathbb{B}^{-1} \left(\ID - \mathbb{A}\right) \traceDe(\ptot - \pinc)
	= -\SL^{-1} \left(\tfrac12\ID - \DL\right) \traceDe\psca
	= \SL^{-1} \left(\DL - \tfrac12\ID\right) \traceDe\psca \\
	&= \DtNe \traceDe\psca
	= \traceNe\psca
	= \traceNe(\ptot - \pinc),
\end{align*}
and thus represents the Neumann trace of the scattered field. Hence,
\begin{align}
	\ptot &= \DLP(\traceDe\ptot) - \SLP(\traceNe\ptot) + \pinc \nonumber \\
	&= \DLP(\traceDi\pint) - \SLP(\vartheta) - \SLP(\traceNe\pinc) + \pinc \qquad \text{in } \Omega^+,
\end{align}
which will be needed for the reconstruction of the exterior field. Thirdly, the unknown surface variable $\sigma$ can be interpreted as a penalty term: this variable will be zero if the regulariser $\REG$ would be exact, i.e., without any numerical errors~\cite{hiptmair2006stabilized}.

\subsection{Regularisation}
\label{sec:regularisation}

The stabilised FEM-BEM algorithm requires choosing a regulariser $\REG$ as defined in Eq.~\eqref{eq:def:reg}. A valid example of a regularisation operator is
\begin{equation}
	\label{eq:regularisation:mh}
	\REG_\text{MH} = \left(\ID - \LB\right)^{-1}
\end{equation}
a Modified Helmholtz (MH) operator~\cite{hiptmair2006stabilized, casenave2014coupled, gatica2014coupling}, where $\LB$ denotes the Laplace-Beltrami operator on the surface~$\Gamma$. The sesquilinear form of its inverse (see Eq.~\eqref{eq:b}) is given by
\begin{equation}
	\mathsf{b}_\text{MH}(\sigma,\tau) = \ipg{\nabla\sigma}{\nabla\tau} + \ipg{\sigma}{\tau} \qquad \text{for } \sigma,\tau \in H^1(\Gamma).
\end{equation}
Other choices of regularisation include a localised inverse Laplace-Beltrami operator~\cite{buffa2005regularized} and single-layer operators~\cite{buffa2006acoustic}. In this study, two different regularisation operators will be introduced.

\subsubsection{Shifted Laplace regularisation}

Let us define
\begin{equation}
	\label{eq:regularisation:sl}
	\REG_\mathrm{SL} = \left(\kappa^2\ID - \LB\right)^{-1}
\end{equation}
a Shifted Laplace (SL) operator, for $\kappa\in\mathbb{R}$ and $\kappa>0$. The associated sesquilinear form is
\begin{equation}
	\mathsf{b}_\mathrm{SL}(\sigma,\tau) = \ipg{\nabla\sigma}{\nabla\tau} + \ipg{\kappa\sigma}{\kappa\tau} \qquad \text{for } \sigma,\tau \in H^1(\Gamma).
\end{equation}
Notice that
\begin{align*}
	\ipg{\theta}{\REG_\mathrm{SL}\theta}
	&= \mathsf{b}_\mathrm{SL}(\REG_\mathrm{SL}\theta,\REG_\mathrm{SL}\theta)
	= \ipg{\nabla(\REG_\mathrm{SL}\theta)}{\nabla(\REG_\mathrm{SL}\theta)} + \ipg{\kappa \REG_\mathrm{SL}\theta}{\kappa \REG_\mathrm{SL}\theta} \\
	&= \normg{\nabla(\REG_\mathrm{SL}\theta)}^2 + \normg{\kappa \REG_\mathrm{SL}\theta}^2
\end{align*}
so that $\ipg{\theta}{\REG_\mathrm{SL}\theta} > 0$ for all $\theta \in H^{-1/2}(\Gamma)\setminus\{0\}$ and the regulariser is a valid choice as stabilisation operator. While $\kappa=1$ reduces to the Modified Helmholtz regulariser, here $\kappa=\kext$ will be chosen.

\subsubsection{OSRC regularisation}

Another idea for regularisation is to use the NtD map. However, no explicit expressions are available for general cases and approximations need to be used. The OSRC-NtD operator is an approximation of the NtD map based on on-surface radiation conditions and is especially accurate at high frequencies. The OSRC-NtD operator is defined by
\begin{equation}
	\label{eq:osrc:ntd}
	\OsrcNtD = \frac1{\imath k} \left(\ID + \frac1{k_\epsilon^2} \LB\right)^{-1/2}
\end{equation}
for $k_\epsilon = k(1+\epsilon\imath )$ with $\epsilon\in\mathbb{R}$ and $\epsilon>0$ (c.f.~\cite{antoine2021introduction}). The damping parameter~$\epsilon$ avoids singularities, the square root operation is approximated by a Padé series expansion, and the remaining inverse operations will be handled on a discrete level by taking a sparse LU factorisation of the surface Helmholtz problem~\cite{darbas2013combining, betcke2017computationally}.

The OSRC-NtD operator is a pseudodifferential operator that maps from $\Hminushalf$ to $\Hplushalf$, c.f.~\cite{antoine2007generalized}.
The definiteness of the operator follows from a spectral analysis. For
\begin{equation}
	-\LB \varphi_j = \mu_j \varphi_j
\end{equation}
the eigendecomposition of the Laplace-Beltrami operator, one has $\mu_0=0$ and $\mu_j > 0$ for $j=1,2,3,\dots$ and $\{\varphi_j\}_{j \in \mathbb{N}}$ an orthonormal basis~\cite{antoine2007generalized, darbas2004thesis}. For any function $\psi\in\Hminushalf$, one can use a series expansion
\begin{equation}
	\psi = \sum_{j=0}^\infty a_j \varphi_j
\end{equation}
with $a_j \in \mathbb{C}$, and the OSRC-NtD operator satisfies
\begin{equation}
	\OsrcNtD \psi = \frac1{\imath k} \sum_{j=0}^\infty \left( 1 - \frac{\mu_j}{k_\epsilon^2} \right)^{-1/2} a_j \varphi_j
\end{equation}
as demonstrated in~\cite{antoine2007generalized, darbas2004thesis}. Since $\varphi_j$ is an orthonormal basis,
\begin{equation}
	\ipg{\psi}{\OsrcNtD \psi} = \frac1{\imath k} \sum_{j=0}^\infty \left( 1 - \frac{\mu_j}{k_\epsilon^2} \right)^{-1/2} |a_j|^2.
\end{equation}
Notice that
\begin{align*}
	\sqrt{1 - \frac{\mu_j}{k_\epsilon^2}}
	&= \frac1{|k_\epsilon|^2} \sqrt{|k_\epsilon|^4 - k^2 (1 - \epsilon^2) \mu_j + (2 k^2 \epsilon \mu_j) \imath}.
\end{align*}
Since $\epsilon>0$ and $\mu_j \geq 0$, also $(2 k^2 \epsilon \mu_j) \geq 0$, so that $\Im\left(\sqrt{1 - \frac{\mu_j}{k_\epsilon^2}}\right) \geq 0$.
Furthermore,
\begin{equation*}
	\Re\left(\ipg{\psi}{\OsrcNtD \psi}\right) = -\frac1{k} \sum_{j=0}^\infty \frac{\Im\left(\sqrt{1 - \frac{\mu_j}{k_\epsilon^2}}\right)}{\left|\sqrt{1 - \frac{\mu_j}{k_\epsilon^2}}\right|^2} |a_j|^2.
\end{equation*}
This analysis concludes that
\begin{equation}
	\Re\left(\ipg{\psi}{\OsrcNtD \psi}\right) \leq 0
\end{equation}
for all $\psi \in \Hminushalf$. Equality holds if and only if $a_j=0$ for all $j>0$, that is, for $\psi = a_0 \varphi_0$ a constant function. Since $\Im\left(\ipg{\psi}{\OsrcNtD \psi}\right) < 0$ for all $\psi \in \Hminushalf \setminus \{0\}$ one can define an arbitrary small complex rotation $0 < \alpha \ll \pi/2$ such that $\Re\left(\ipg{\psi}{e^{-\imath\alpha}\OsrcNtD \psi}\right) < 0$ for all $\psi \in \Hminushalf \setminus \{0\}$. In practice, this perturbation is not necessary since the discretised variable $\psi$ is not expected to be exactly constant in finite-precision arithmetic. Therefore, we will consider the regulariser
\begin{equation}
	\label{eq:regularisation:ntd}
	\REGntd = -\OsrcNtD = -\frac1{\imath k} \left(\ID + \frac1{k_\epsilon^2} \LB\right)^{-1/2}
\end{equation}
the negative OSRC-NtD operator~\eqref{eq:osrc:ntd}. Its inverse is given by the negative OSRC-DtN operator
\begin{equation}
	\label{eq:regularisation:dtn}
	\left(\REGntd\right)^{-1} = -\OsrcDtN = -\imath k \left(\ID + \frac1{k_\epsilon^2} \LB\right)^{1/2},
\end{equation}
whose weak formulation can directly be used for the sesquilinear form $\mathsf{b}_\mathrm{OSRC}$.

\subsection{Discrete system}

The variational problems of the FEM-BEM formulations will be discretised with a Galerkin method on a tetrahedral mesh in $\Omega^-$ and a triangular mesh on $\Gamma$ that match at the material interface (cf.~\cite{logg2012automated, smigaj2015solving}). The standard FEM-BEM algorithm (also known as Johnson-Nédélec coupling~\cite{johnson1980coupling, sayas2009validity}) uses the DtN map~\eqref{eq:dtn:standard} and results in the linear system
\begin{equation}
	\label{eq:system:standard}
	\begin{bmatrix}
		\mathcal{F} & -\frac\rhoext\rhoint \ID \\
		\tfrac12\ID - \DL & \SL \end{bmatrix}
	\begin{bmatrix} \pint \\ \theta \end{bmatrix}
	= \begin{bmatrix}
		0 \\
		\traceDe\pinc
	\end{bmatrix}
\end{equation}
where $\mathcal{F}$ denotes the discrete sesquilinear form $\mathsf{a}(\pint,q)$ for the FEM in the interior domain. The discrete variational problem for the symmetric FEM-BEM~\eqref{eq:fembem:symmetric} is given by the system matrix
\begin{equation}
	\label{eq:system:symmetric}
	\begin{bmatrix}
		\mathcal{F} + \frac\rhoext\rhoint \HS & \frac\rhoext\rhoint (\AD - \tfrac12\ID) \\
		\tfrac12\ID - \DL & \SL \end{bmatrix}
	\begin{bmatrix} \pint \\ \theta \end{bmatrix}
	= \begin{bmatrix}
		\frac\rhoext\rhoint \HS \traceDe\pinc + \frac\rhoext\rhoint \traceNe\pinc \\
		(\tfrac12\ID - \DL) \traceDe\pinc
	\end{bmatrix}
\end{equation}
The discrete variational problem for the stabilised FEM-BEM~\eqref{eq:fembem:stable} is given by
\begin{align}
	\label{eq:system:stable}
	&\begin{bmatrix}
		\mathcal{F} + \frac\rhoext\rhoint \left(\HS + \imath\nu(\tfrac12\ID - \DL)\right) & \frac\rhoext\rhoint \left(\AD - \tfrac12\ID + \imath\nu\SL\right) & 0 \\
		\tfrac12\ID - \DL & \SL & \imath\eta\ID \\
		-\HS & -(\tfrac12\ID + \AD) & \mathcal{S}
	\end{bmatrix}
	\begin{bmatrix} \pint \\ \theta \\ \Sigma \end{bmatrix} \nonumber \\
	&\qquad = \begin{bmatrix}
		\frac\rhoext\rhoint \left(\HS + \imath\nu(\tfrac12\ID - \DL)\right) \traceDe\pinc + \frac\rhoext\rhoint \traceNe\pinc \\
		(\tfrac12\ID - \DL) \traceDe\pinc \\
		-\HS \traceDe\pinc
	\end{bmatrix}
\end{align}
where $\mathcal{S}$ denotes the discrete sesquilinear form $\mathsf{b}(\sigma,\tau)$ for the inverse regulariser.

\subsubsection{Alternative formulations}

Notice that the stabilised FEM-BEM formulation~\eqref{eq:system:stable} can also be written as
\begin{align}
	\label{eq:system:stable:alternative:nu}
	&\begin{bmatrix}
		\mathcal{F} + \frac\rhoext\rhoint \HS & \frac\rhoext\rhoint (\AD - \tfrac12\ID) & \frac\rhoext\rhoint \nu\eta\ID \\
		\tfrac12\ID - \DL & \SL & \imath\eta\ID \\
		-\HS & -(\tfrac12\ID + \AD) & \mathcal{S}
	\end{bmatrix}
	\begin{bmatrix} \pint \\ \theta \\ \Sigma \end{bmatrix}
	= \begin{bmatrix}
		\frac\rhoext\rhoint \HS \traceDe\pinc + \frac\rhoext\rhoint \traceNe\pinc \\
		(\tfrac12\ID - \DL) \traceDe\pinc \\
		-\HS \traceDe\pinc
	\end{bmatrix}
\end{align}
which has the symmetric FEM-BEM formulation in the upper left blocks. This alternative formulation also emphasises the contribution of the stabilisation term $\Sigma$ as a penalty.

Finally, for $\sigma=\REG\hat\sigma$, the stabilised formulation~\eqref{eq:system:stable} can be written as
\begin{align}
	\label{eq:system:stable:alternative:reg}
	&\begin{bmatrix}
		\mathcal{F} + \frac\rhoext\rhoint \HS & \frac\rhoext\rhoint (\AD - \tfrac12\ID) & \frac\rhoext\rhoint \nu\eta \REG \\
		\tfrac12\ID - \DL & \SL & \imath\eta \REG \\
		-\HS & -(\tfrac12\ID + \AD) & \ID
	\end{bmatrix}
	\begin{bmatrix} \pint \\ \theta \\ \hat\Sigma \end{bmatrix}
	= \begin{bmatrix}
		\frac\rhoext\rhoint \HS \traceDe\pinc + \frac\rhoext\rhoint \traceNe\pinc \\
		(\tfrac12\ID - \DL) \traceDe\pinc \\
		-\HS \traceDe\pinc
	\end{bmatrix}
\end{align}
which is more convenient in case the regulariser is explicitly defined.

\subsubsection{Restriction operators}

Notice that on the first row of the FEM-BEM formulations, both FEM and BEM terms are present. The BEM operators only act on a subspace of $H^1(\Omega^-)$, namely $\Hplushalf$ where $\Gamma = \partial\Omega^-$. For this purpose, let us define the operator $\TR$ that discretises the Dirichlet trace operation $\traceDi: H^1(\Omega^-) \to \Hplushalf$. This linear operator restricts the FEM's degrees of freedom to the corresponding degrees of freedom in the BEM, and can be represented by a sparse matrix. Then, system~\eqref{eq:system:stable} reads
\begin{align}
	\label{eq:system:stable:restrictions}
	&\begin{bmatrix}
		\mathcal{F} + \frac\rhoext\rhoint \TR^T \left(\HS + \imath\nu(\tfrac12\ID - \DL)\right) \TR & \frac\rhoext\rhoint \TR^T \left(\AD - \tfrac12\ID + \imath\nu\SL\right) & 0 \\
		\left(\tfrac12\ID - \DL\right) \TR & \SL & \imath\eta\ID \\
		-\HS \TR & -(\tfrac12\ID + \AD) & \mathcal{S}
	\end{bmatrix}
	\begin{bmatrix} \pint \\ \theta \\ \Sigma \end{bmatrix} \nonumber \\
	&\quad= \begin{bmatrix}
		\frac\rhoext\rhoint \TR^T \left(\HS + \imath\nu(\tfrac12\ID - \DL)\right) \traceDe\pinc + \frac\rhoext\rhoint \TR^T \traceNe\pinc \\
		(\tfrac12\ID - \DL) \traceDe\pinc \\
		-\HS \traceDe\pinc
	\end{bmatrix}.
\end{align}
In the following, the operators $\TR$ will be dropped where possible for brevity.

\subsubsection{Function spaces}

Table~\ref{table:spaces} summarises the function spaces of the operators in the variational formulation, where the test space is the dual of the range space. The Galerkin discretisation takes basis and test functions from finite-dimensional subspaces of the domain and test spaces, respectively. Specifically, the space $H^1(\Omega^-)$ will be approximated by continuous piecewise linear (P1) functions on a tetrahedral mesh that equal one in a single node and zero in all other nodes. On the surface, a triangular mesh that matches the tetrahedral mesh will be used, with node-based continuous piecewise linear (P1) functions for the spaces $H^{1}(\Gamma)$ and $\Hplushalf$. The space $\Hminushalf$ can be discretised with triangle-based piecewise constant (P0) functions that equal one in a triangle and zero in all other triangles. Alternatively, P1 functions can also be used for the space $\Hminushalf$.

\begin{table}[!ht]
	\centering
	\caption{Function spaces of the stabilised FEM-BEM system~\eqref{eq:system:stable}.}
	\label{table:spaces}
	\begin{tabular}{llll}
		\hline\hline
		variable & $H_\text{domain}$ & $H_\text{range}$ & $H_\text{test}$ \\
		\hline
		row/column 1 (FEM) & $H^1(\Omega^-)$ & $H^0(\Omega^-)$ & $H^1(\Omega^-)$ \\
		row/column 1 (BEM) & $H^{1/2}(\Gamma)$ & $H^{-1/2}(\Gamma)$ & $H^{1/2}(\Gamma)$ \\
		row/column 2 & $H^{-1/2}(\Gamma)$ & $H^{1/2}(\Gamma)$ & $H^{-1/2}(\Gamma)$ \\
		row/column 3 & $H^{1/2}(\Gamma)$ & $H^{-1/2}(\Gamma)$ & $H^{1/2}(\Gamma)$ \\
		row/column 3 \eqref{eq:system:stable:alternative:reg} & $H^{-1/2}(\Gamma)$ & $H^{-1/2}(\Gamma)$ & $H^{1/2}(\Gamma)$ \\
		\hline\hline
	\end{tabular}
\end{table}

\subsection{Preconditioning}
\label{sec:preconditioning}

The weak formulation of the FEM-BEM variational problem needs to be solved with a linear solver. Direct solvers have a high computational complexity and are difficult to combine with matrix compression techniques. Hence, they become intractable for large-scale simulations. Alternatively, iterative solvers only need an implementation of a matrix-vector multiplication and have a low computational complexity. Here, the GMRES algorithm will be used because the system matrix is indefinite. The convergence of the GMRES depends on the spectrum of the system matrix and can be improved significantly with preconditioning~\cite{antoine2021introduction}. Below, different preconditioning strategies will be considered.

\subsubsection{Permutation}

One of the easiest preconditioning strategies is to consider permutations of (blocks of) rows in the system matrix. While this does not reduce the condition number, the spectrum changes and the convergence of linear solvers can improve considerably. Importantly, matrix permutations do not affect the computational costs of a matrix-vector multiplication. Specifically, swapping the last two block rows in the stabilised formulation~\eqref{eq:system:stable:alternative:nu} yields
\begin{align}
	\label{eq:system:stable:alternative:nu:permuted}
	&\begin{bmatrix}
		\mathcal{F} + \frac\rhoext\rhoint \HS & \frac\rhoext\rhoint (\AD - \tfrac12\ID) & \frac\rhoext\rhoint \nu\eta\ID \\
		\HS & \tfrac12\ID + \AD & -\mathcal{S} \\
		\tfrac12\ID - \DL & \SL & \imath\eta\ID
	\end{bmatrix}
	\begin{bmatrix} \pint \\ \theta \\ \Sigma \end{bmatrix}
	= \begin{bmatrix}
		\frac\rhoext\rhoint \HS \traceDe\pinc + \frac\rhoext\rhoint \traceNe\pinc \\
		\HS \traceDe\pinc \\
		(\tfrac12\ID - \DL) \traceDe\pinc
	\end{bmatrix}
\end{align}
where the sign of the second row was changed as well. This formulation is expected to improve the conditioning since the main diagonal includes the compact operators $\tfrac12\ID+\AD$ and $\ID$ instead of the single-layer operator~$\SL$ and the regularisation term~$\mathcal{S}$. In other words, the domain and range spaces are the same for the second and third row/column. Notice that a P1-P1 discretisation is preferred over a P0-P1 discretisation to avoid rectangular blocks on the main diagonal.
The variational form for the FEM (denoted by $\mathcal{F}$) has to remain on the main diagonal because any permutation will leave a pure BEM operator in the top left block. These BEM operators are restricted to the surface unknowns of $\pint$ and have rows and columns with only zeros.

\subsubsection{Combined FEM and BEM preconditioning}
\label{sec:prec:fembem}

The top left block of the stabilised FEM-BEM formulations is a linear combination of a FEM and BEM operator, where the BEM operator only acts on the degrees of freedom on the material interface. Hence, two preconditioning strategies can be distinguished: 1) a single preconditioner for the entire block, and 2) a separate preconditioner for the surface and interior nodes of the volume mesh. The second choice can be written as
\begin{equation*}
	\begin{bmatrix} \TRc^T & \TR^T \end{bmatrix} \begin{bmatrix} M_\mathrm{FEM} & 0 \\ 0 & M_\mathrm{BEM} \end{bmatrix} \begin{bmatrix} \TRc \\ \TR \end{bmatrix} \left(\mathcal{F} + \frac\rhoext\rhoint \TR^T \HS \TR\right) \pint
\end{equation*}
where $\TR$ denotes the restriction to the surface nodes and $\TRc$ to the interior nodes of the volumetric mesh. Notice that this does not split the FEM operator entirely from the BEM operator: at the boundary a linear combination remains present. Hence, it is a design choice to use either the FEM-based or the BEM-based preconditioner for this combined term.

The difficulty of a single preconditioner for the combined FEM and BEM block is that most preconditioners designed for FEM are not adequate for the BEM. For example, almost all FEM preconditioners require the matrix elements to be available explicitly, which is not the case for compression techniques (e.g., fast multipole or hierarchical matrices). Also, FEM preconditioners typically take the sparsity of the matrix into account while the BEM matrix is dense. One could design preconditioners based on sparsified BEM matrices~\cite{liang2019coupled, gao2016non, bruckner20123d, gaul2008coupling}, but this avenue will not be pursued here.

Generally speaking, the number of GMRES iterations that are required to solve the symmetric formulation~\eqref{eq:system:symmetric} increases quickly with model complexity and optimal Schwarz preconditioners can be designed for the full block version~\cite{stephan2018coupling}. However, these preconditioners become prohibitively expensive for large-scale simulations due to the inversion of operator blocks or decompositions of the discrete space~\cite{hahne1995fast, funken1996hierarchical, heuer1998preconditioners, heuer1999preconditioned, kuhn2002symmetric, harbrecht2003multiscale, vouvakis2007domain, funken2009fast, feischl2017optimal}. Similarly, preconditioning based on domain decomposition~\cite{langer1994parallel} increases the implementation effort and restricts the flexibility of the FEM-BEM coupling algorithm.

Here, operator preconditioning will be used for the BEM operators and standard incomplete factorisations for the FEM operator, either for the interior nodes or for the entire operator.

\subsubsection{BEM preconditioning}
\label{sec:prec:bem}

Operator preconditioning is a successful strategy for boundary integral operators where the preconditioner is defined as a continuous operator and discretised separately from the model. Advantages of operator preconditioning include the availability of a mathematical foundation~\cite{hiptmair2006operator, kirby2010functional} and easy combination with matrix compression techniques~\cite{peeters2010embedding, niino2012preconditioning, darbas2013combining}. Furthermore, few numerical parameters need to be set and its effectiveness has been shown for large-scale simulations of different wave phenomena (cf.~\cite{andriulli2008multiplicative, dobbelaere2015calderon, chaillat2017fast, haqshenas2021fast, wout2021benchmarking}).

\paragraph{Mass preconditioning}
While the strong formulation of a boundary integral operator maps from $H_\text{domain}$ to $H_\text{range}$, the weak formulation actually maps from $H_\text{domain}$ to $H_\text{test}$. The idea of mass preconditioning is to take the inverse of an identity operator that maps from $H_\text{range}$ to $H_\text{test}$ as preconditioner~\cite{hiptmair2006operator, kirby2010functional}. The advantage is that conditioning is improved while the sparse LU factorisation of the mass matrix yields little computational overhead. Here, mass preconditioning will only be used for P1-P1 discretisation. In the case of P0-P1 discretisation, the mass matrix is rectangular. Furthermore, no P0 functions on dual meshes will be considered to avoid computational overhead in the quadrature scheme~\cite{betcke2020product}.

\paragraph{OSRC preconditioning}
The hypersingular and single-layer boundary integral operators can be preconditioned with operators of opposite order~\cite{steinbach1998construction}. Among the most effective choices of preconditioners are the OSRC-approximated NtD and DtN maps, denoted by $\OsrcNtD$ and $\OsrcDtN$ and defined by Eqns.~\eqref{eq:regularisation:ntd} and~\eqref{eq:regularisation:dtn}, respectively. They are efficient because of their sparsity and are accurate at high frequencies~\cite{antoine2021introduction, wout2021pmchwt}. An example of OSRC preconditioning for the FEM-BEM system (Eq.~\eqref{eq:system:stable} with $\nu=0$) is
\begin{align*}
	&\begin{bmatrix} \TRc^T \TRc + \TR^T \OsrcNtD \TR & 0 & 0 \\ 0 & \OsrcDtN & 0 \\ 0 & 0 & \OsrcNtD \end{bmatrix}
	\begin{bmatrix}
		\mathcal{F} + \frac\rhoext\rhoint \TR^T \HS \TR & \frac\rhoext\rhoint \TR^T (\AD - \tfrac12\ID) & 0 \\
		\left(\tfrac12\ID - \DL\right) \TR & \SL & \imath\eta\ID \\
		-\HS \TR & -(\tfrac12\ID + \AD) & \mathcal{S}
	\end{bmatrix}
	\begin{bmatrix} \pint \\ \theta \\ \Sigma \end{bmatrix}.
\end{align*}
The OSRC operators require a P1-P1 discretisation. Also, while the stabilisation term $\mathcal{S}$ is not a hypersingular boundary integral, the OSRC-NtD operator is still of opposite order, that is, $\OsrcNtD \mathcal{S}: \Hplushalf \to \Hplushalf$.

\subsubsection{FEM preconditioning}

The operator $\mathcal{F}$ for the FEM contribution is sparse and its dimensions are larger than for the BEM operators. Furthermore, the number of degrees of freedom in the FEM scales cubic with frequency while the BEM scales quadratic. Hence, a complete factorisation of $\mathcal{F}$ becomes intractable for large-scale simulations and sparse preconditioners need to be employed. One of the most effective approaches is incomplete LU (ILU) preconditioning, which sets up a sparse approximation of the inverse matrix. Fill-in elements are dropped based on a predefined sparsity pattern or a given threshold. Here, the SuperLU algorithm will be used~\cite{li2011supernodal}. Other FEM preconditioners (cf.~\cite{logg2012automated}) can easily be used as well but will not be considered for brevity.

\section{Results}
\label{sec:results}

This section showcases numerical results and investigates the influence of model choices on the computational performance of the stabilised FEM-BEM algorithm.

\subsection{Benchmark settings}
\label{sec:results:settings}

\subsubsection{Computational platform}

The coupled FEM-BEM algorithm was implemented in BEMPP~\cite{smigaj2015solving, betcke2021bempp} (version 0.2 of bempp-cl) and FEniCS~\cite{logg2012automated} (version 2019.1). The linear algebra was performed with SciPy~\cite{virtanen2020scipy} (version 1.6) and includes the GMRES algorithm~\cite{saad1986gmres} and SuperLU factorisation~\cite{li2011supernodal}. Multithreading was enabled through the Numba implementation~\cite{lam2015numba} of BEMPP.

The numerical parameters are the default values of the respective algorithms and libraries. The GMRES solver has a termination criterion of $10^{-5}$ and no restart. The ILU factorisation has a drop tolerance of $10^{-4}$. The OSRC operators use a Padé expansion of size two, a branch cut of $\pi/3$, and a damping parameter of $\epsilon = 0.4 (\kext L)^{-\frac23}$ where $L$ is the characteristic length of the scatterer~$\Omega^-$ (cf.~\cite{antoine2021introduction}).

\subsubsection{Geometry}

A unit cube will be used as geometry since the resonance modes are known explicitly for this structure. That is,
\begin{equation}
	\label{eq:modes:box}
	k_m = \pi \sqrt{m_x^2 + m_y^2 + m_z^2} \quad \text{ for } m_x, m_y, m_z \in \{1,2,3,\dots\}
\end{equation}
are the wavenumbers at which resonances occur. Let us consider materials that satisfy a constant density ($\rhoint=\rhoext$) and an interior wavenumber $k(\mathbf{x}) = \kext n(\mathbf{x})$ for a given exterior wavenumber $\kext$ and refractivity $n(\mathbf{x})$ for $\mathbf{x}\in\Omega^-$. The refractivity is chosen as
\begin{equation} \label{eq:refractivity}
	n(\mathbf{x}) = \frac{1 - \tfrac12 e^{-\norm{\mathbf{x} - \tfrac12}_\infty^2}}{1 - \tfrac12 e^{-1/4}}.
\end{equation}
The incident field is a plane wave with direction vector $(1/\sqrt{5},2/\sqrt{5},0)$.

The interior of the cube is meshed with tetrahedral elements. The triangular mesh on the surface of the cube is extracted from the tetrahedral elements. The grid has 2744 degrees of freedom for the FEM space, 2028 for the P0 space, and 1016 for the P1 space in BEM. The maximum cell size is 0.133 so that at least four elements per wavelength are present at a wavenumber of 12.

\FloatBarrier
\subsection{Acoustic field at a resonating cube}

The salient feature of the stabilised FEM-BEM algorithm is its stability at resonance frequencies. Let us validate the stability at a unit cube, for which the resonance modes are given by Eq.~\eqref{eq:modes:box}. A direct solver was used for this benchmark to avoid numerical errors of iterative linear solvers.

\begin{figure}[!ht]
	\centering
	\includegraphics[width=.9\columnwidth]{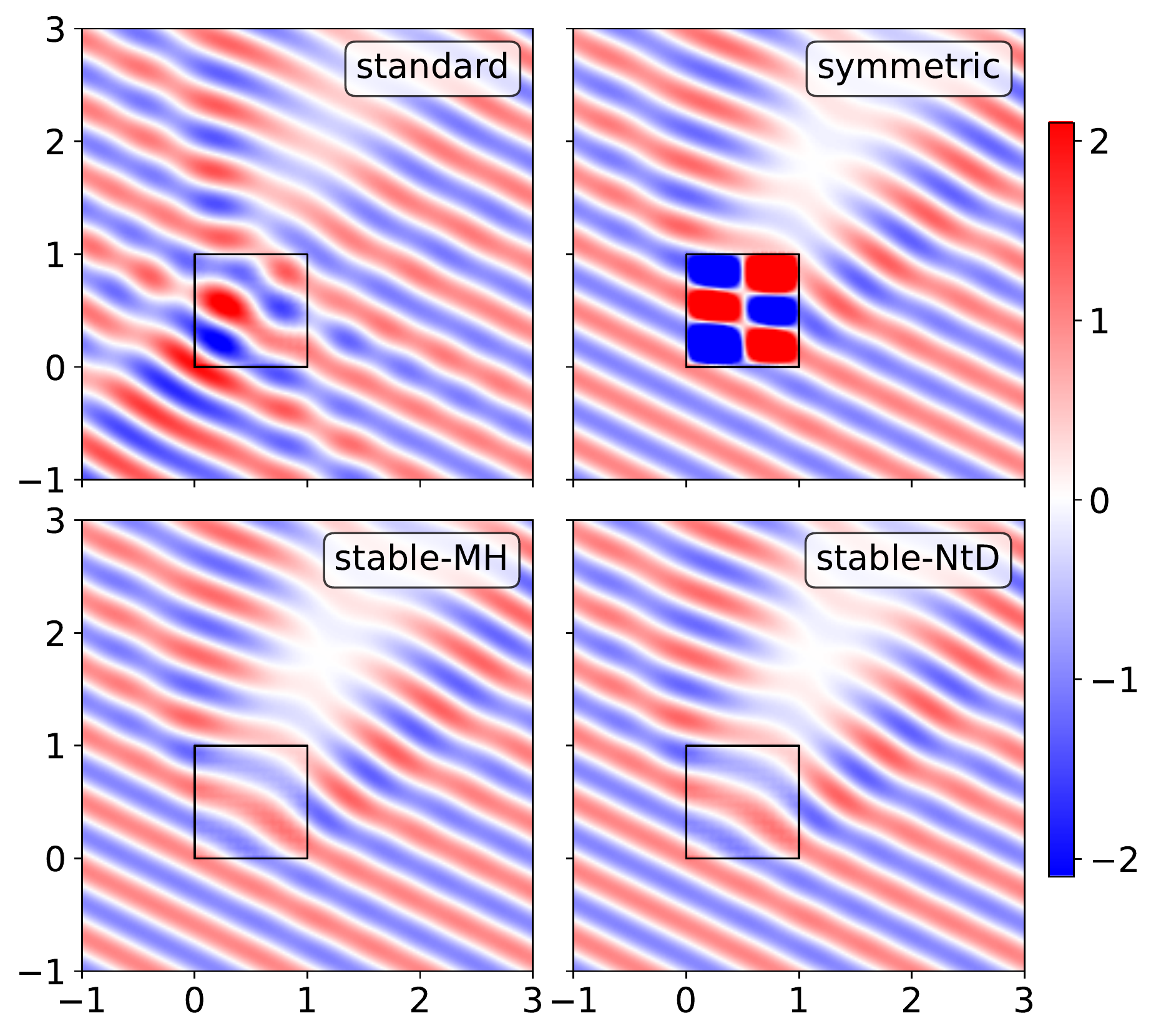}
	\caption{The real part of the acoustic field for different FEM-BEM coupling techniques at $\kext=11.7519$, on the plane $z=0.5$. The stabilisation parameters are $\eta=1$ and $\nu=1$ and a direct solver is used to solve the weak formulation with P0-P1 elements. The formulations are: standard~\eqref{eq:system:standard}, symmetric~\eqref{eq:fembem:symmetric}, and stabilised~\eqref{eq:fembem:stable} with Modified Helmholtz~\eqref{eq:regularisation:mh} and OSRC Neumann-to-Dirichlet~\eqref{eq:regularisation:ntd} regularisation.}
	\label{fig:field:cube}
\end{figure}

Figure~\ref{fig:field:cube} depicts the acoustic field calculated with different FEM-BEM formulations. The computational results show the impact of resonance frequencies on the stability of FEM-BEM coupling techniques: the acoustic fields from the standard and symmetric FEM-BEM systems are inaccurate. Notice that the exterior wavenumber of 11.7519 is close to the resonance mode of
\begin{equation*}
	k_{(1,2,3)} = \pi \sqrt{1^2 + 2^2 + 3^2} = 11.7548.
\end{equation*}
The block structure of two and three peaks is indeed visible in the field of the symmetric formulation. Differently, the acoustic fields simulated by the stabilised formulations are accurate.

\FloatBarrier
\subsection{Convergence of the linear solver at a resonating cube}

\begin{figure}[!ht]
	\centering
	\begin{subfigure}[b]{\columnwidth}
		\includegraphics[width=\columnwidth]{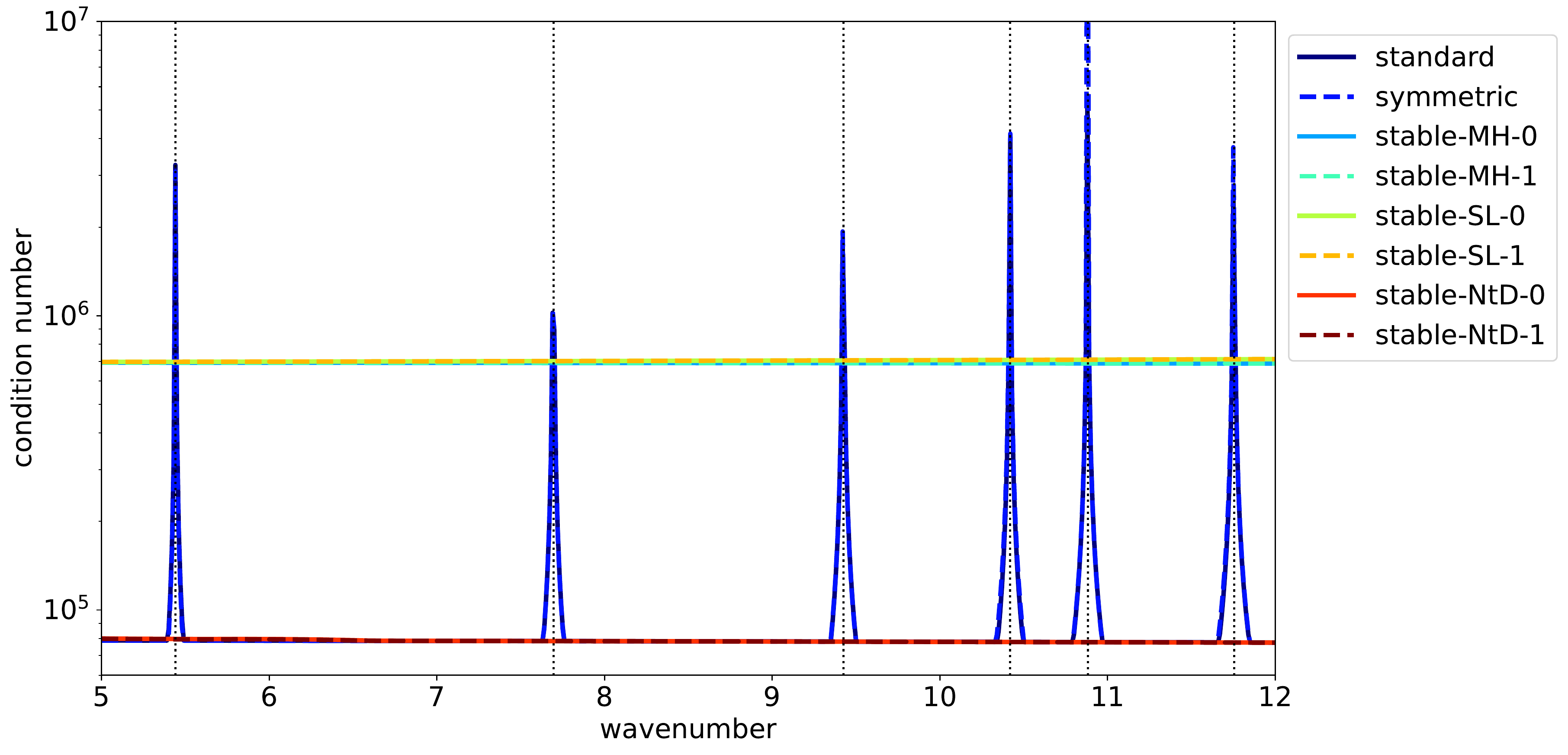}
		\caption{Condition number of the system matrix.}
	\end{subfigure}
	\begin{subfigure}[b]{\columnwidth}
		\includegraphics[width=\columnwidth]{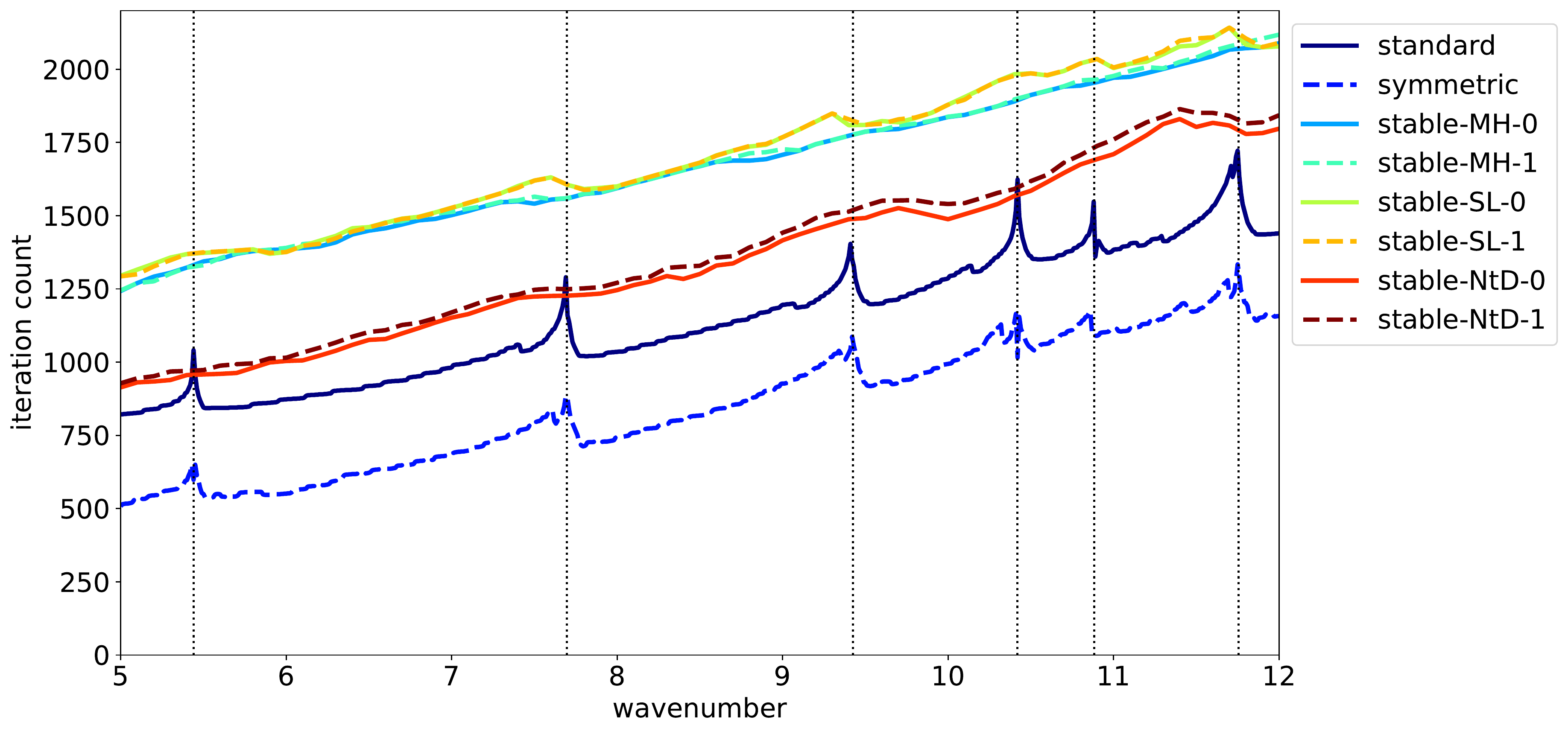}
		\caption{Number of GMRES iterations.}
	\end{subfigure}
	\caption{The conditioning of different FEM-BEM coupling techniques with respect to the exterior wavenumber. The legend specifies the value of $\nu$ and the other stabilisation parameter is $\eta=1$. The BEM uses a P0-P1 discretisation, without preconditioning. The vertical lines depict the analytical resonance frequencies of a unit cube.}
	\label{fig:cube}
\end{figure}

The problem with instability at resonance frequencies is twofold: the simulated field can be inaccurate and the system matrix can be ill conditioned. The poor conditioning deteriorates the convergence of the iterative linear solver. Figure~\ref{fig:cube} shows the condition number and GMRES iteration count for a parameter sweep over the exterior wavenumbers, all on the same mesh. The results show spikes at the resonance frequencies of the cube for the standard and symmetric FEM-BEM algorithm, whereas the stabilised formulations remain well-conditioned at all frequencies.

\subsection{The influence of the model choices on the computational performance}
\label{sec:results:parameters}

The computational results presented above demonstrate the effectiveness of the stabilised FEM-BEM algorithm in suppressing spurious solutions at resonance frequencies. This robustness comes at a price of higher computational costs. Considerably more boundary integral operators need to be assembled and a larger system needs to be solved. The design of the stabilised FEM-BEM algorithm includes various model choices that influence the computational performance. Following is a list of the most influential design choices.
\begin{itemize}
	\item The FEM-BEM formulation: standard~\eqref{eq:system:standard}, symmetric~\eqref{eq:fembem:symmetric}, and stabilised~\eqref{eq:fembem:stable} with Modified Helmholtz (MH), shifted Laplacian (SL) or Neumann-to-Dirichlet (NtD) regularisation.
	\item The choice of stabilisation parameters $\eta$ and $\nu$.
	\item The choice of discretisation based on either P0 or P1 test and basis functions for the space $\Hminushalf$.
	\item The choice of the alternative formulations~\eqref{eq:system:stable:alternative:nu} and~\eqref{eq:system:stable:alternative:reg}.
	\item The choice of the preconditioner (see Section~\ref{sec:preconditioning}).
\end{itemize}
The effect of the design choices on the convergence of GMRES will be benchmarked on a unit cube with constant density and a refractivity as in Eq.~\eqref{eq:refractivity}. The frequency sweep is performed on the same mesh as specified above. In the following, the condition number will not be reported anymore for brevity. Moreover, the condition number is a poor estimate of actual GMRES convergence~\cite{antoine2021introduction}, especially for indefinite systems such as the FEM-BEM formulation.

\begin{figure}[!ht]
	\centering
	\includegraphics[width=\columnwidth]{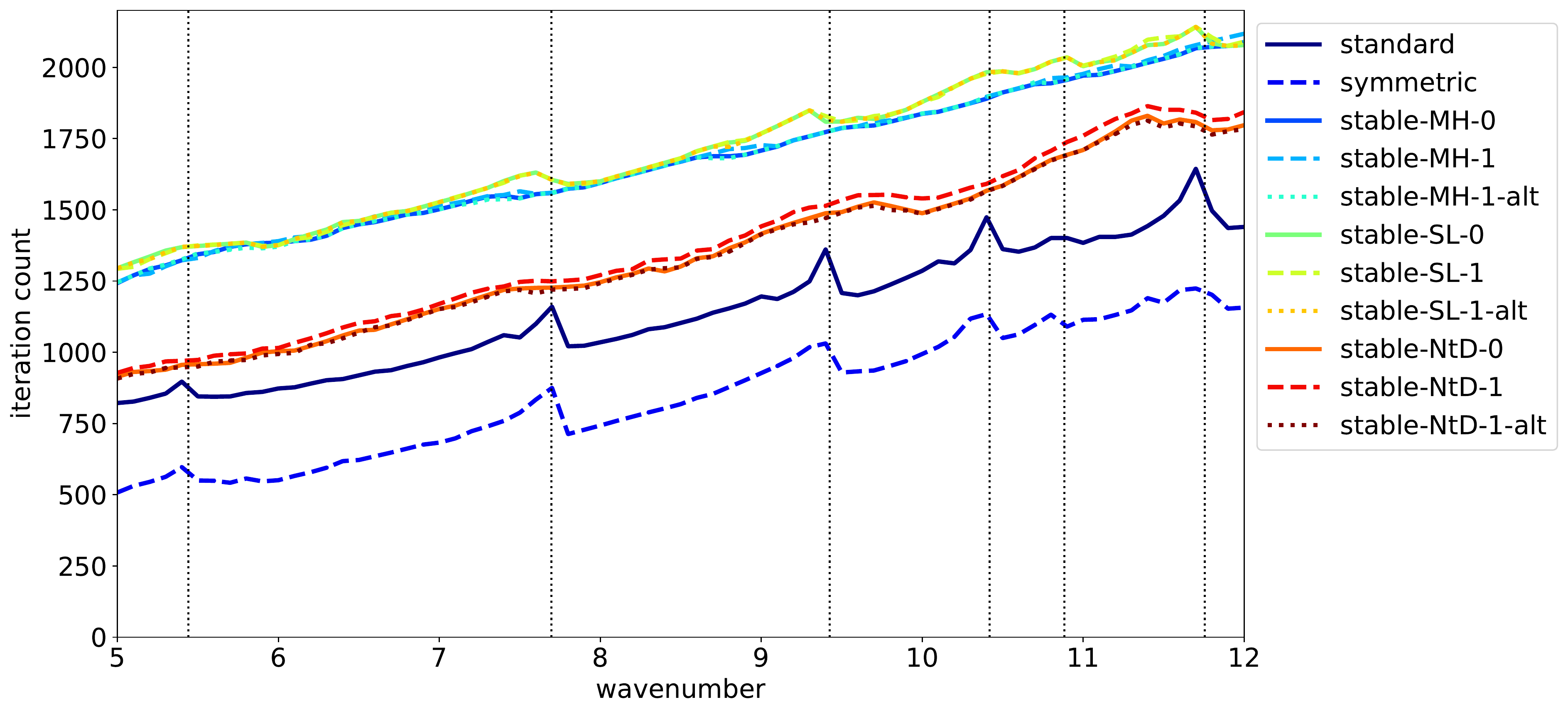}
	\caption{The influence of the stabilisation parameter $\nu$ on the performance of the different FEM-BEM coupling techniques. The legend specifies the value of $\nu$ and \emph{alt} means the alternative formulation~\eqref{eq:system:stable:alternative:nu}. The other stabilisation parameter is $\eta=1$ and a P0-P1 discretisation is used without preconditioning.}
	\label{fig:cube:nu}
\end{figure}

\FloatBarrier
\subsubsection{The stabilisation parameters}

The choice of stabilisation parameters $\eta$ and $\nu$ in the stabilisation operator~\eqref{eq:tracetransformation} is an open question, with a value of $\eta=1$ a good guideline~\cite{meury2007stable}. The parameter $\nu$ can be set to either zero or one. In other words, stabilisation can be included or excluded from the first equation in the variational formulation~\eqref{eq:fembem:stable}. Another model choice that impacts the stabilisation algorithm is the option to use the alternative formulation~\eqref{eq:system:stable:alternative:nu} in which the stabilisation term is present at a different location on the first row.

Figure~\ref{fig:cube:nu} presents the performance of the three cases: 1) $\nu=0$; 2) $\nu=1$ for Eq.~\eqref{eq:system:stable}; and 3) $\nu=1$ for Eq.~\eqref{eq:system:stable:alternative:nu}. Only minor differences are observed between the different settings. Hence, $\nu=0$ will be chosen in the following for efficiency reasons: it requires the smallest number of boundary integral operators while still stabilising the FEM-BEM algorithm.

\FloatBarrier
\subsubsection{The discrete space}

\begin{figure}[!ht]
	\centering
	\includegraphics[width=\columnwidth]{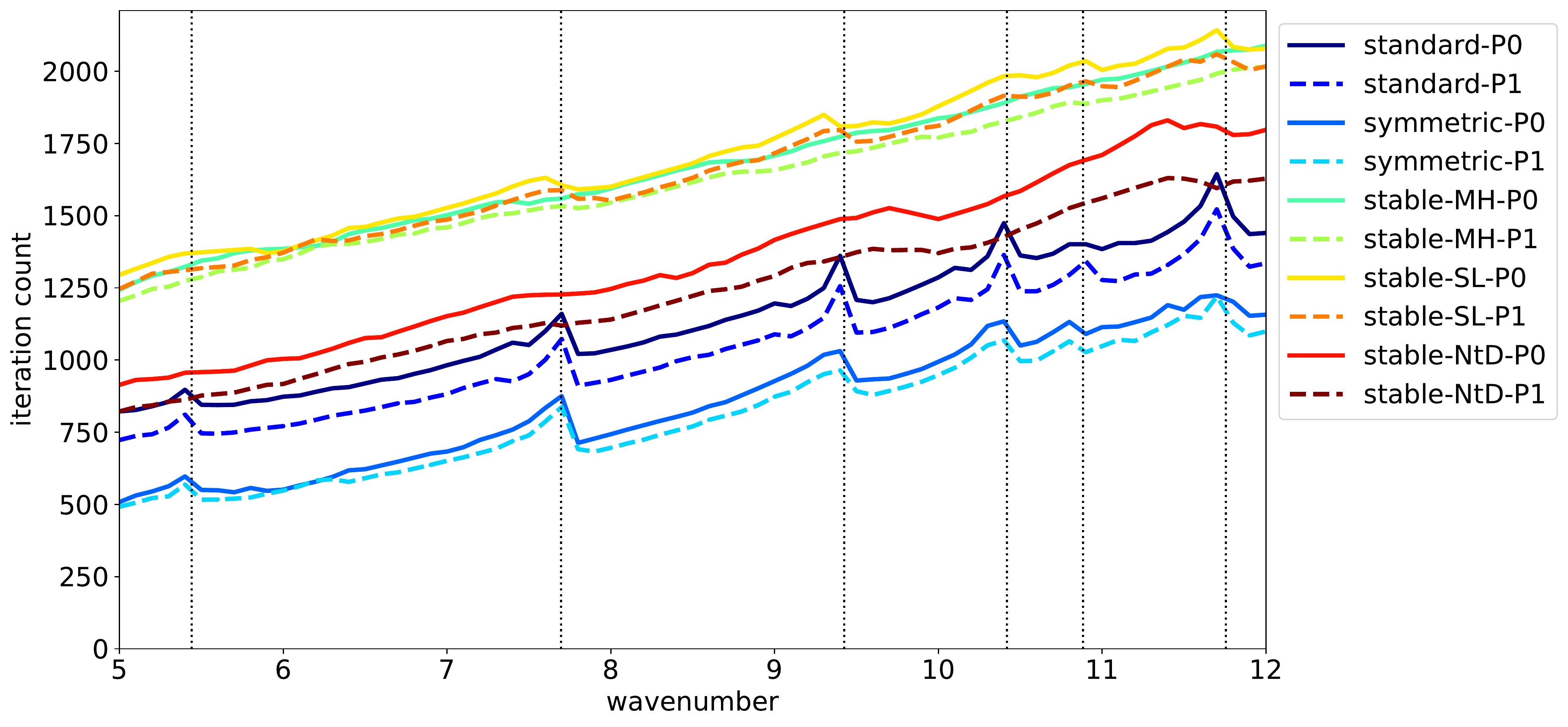}
	\caption{The influence of the discrete space on the performance of the different FEM-BEM coupling techniques. The stabilisation parameters are $\eta=1$ and $\nu=0$, and no preconditioning is used.}
	\label{fig:cube:space}
\end{figure}

The space $\Hminushalf$ can be discretised with either discontinuous P0 elements or continuous P1 elements. The results in Figure~\ref{fig:cube:space} show that, in general, the P1-P1 discretisation of the BEM operators yields a faster convergence of GMRES than the P0-P1 discretisation of the spaces $\Hminushalf$ and $\Hplushalf$. Furthermore, the triangle-based P0 space has approximately twice as many degrees of freedom than the node-based P1 space. Finally, P1 elements allow for operator preconditioning. For these reasons, only P1-P1 discretisations will be considered in the following.

The previous benchmarks showed that the shifted Laplacian regularisation has a very similar performance as the modified Helmholtz regularisation. For clarity of presentation, the SL regulariser will not be shown anymore in the following.

\FloatBarrier
\subsubsection{Operator preconditioning}

\begin{figure}[!ht]
	\centering
	\includegraphics[width=\columnwidth]{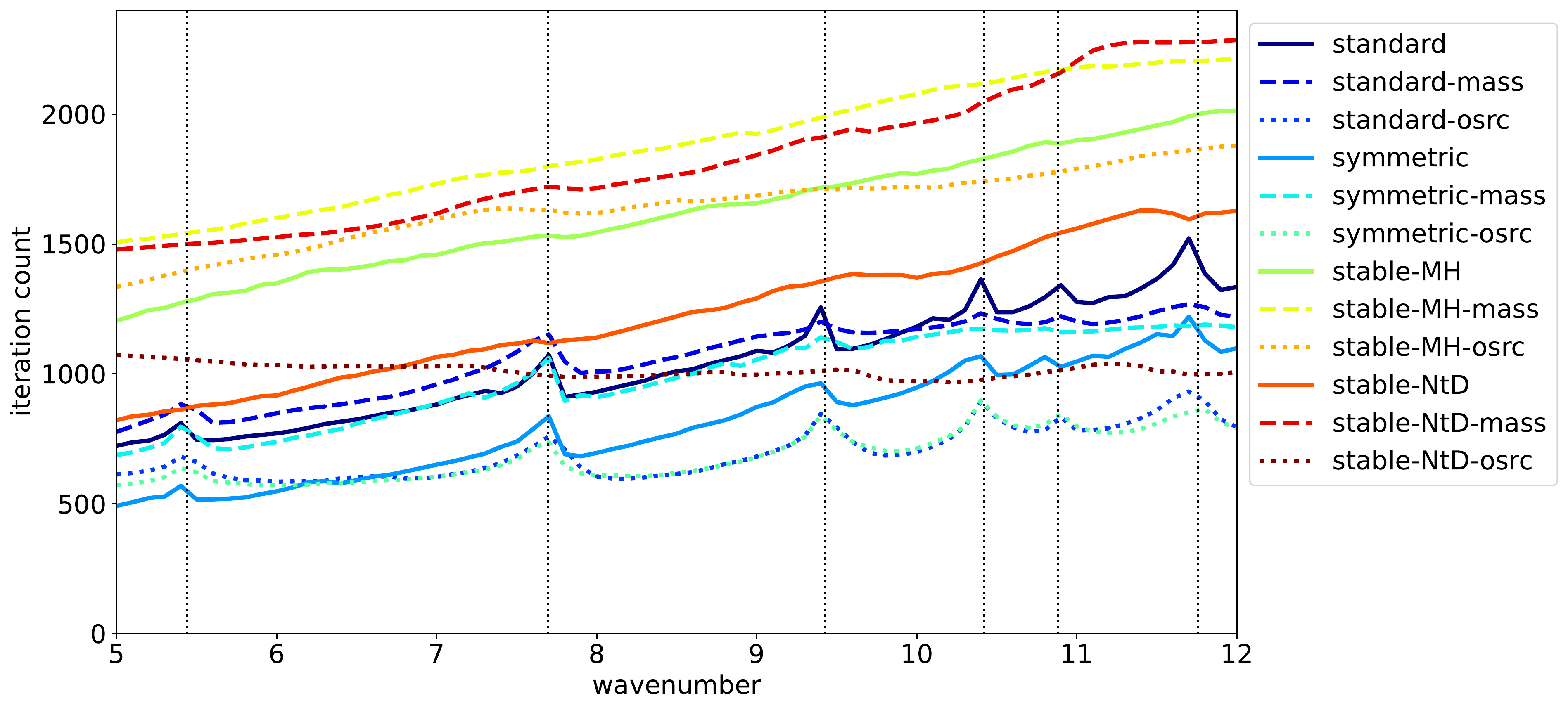}
	\caption{The influence of mass and OSRC preconditioning on the performance of the different FEM-BEM coupling techniques. The stabilisation parameters are $\eta=1$ and $\nu=0$, and P1 elements are used.}
	\label{fig:cube:prec:operator}
\end{figure}

As explained in Section~\ref{sec:prec:bem}, operator preconditioning promises to improve the GMRES convergence. The results in Figure~\ref{fig:cube:prec:operator} show that operator preconditioning can improve but also deteriorate the convergence of GMRES. This behaviour can be attributed to the fact that only the BEM operators are preconditioned, which can yield an imbalance in the preconditioned residual compared to the FEM block. In general, the OSRC preconditioner is effective at reducing the iteration count, especially at higher frequencies where the system is more difficult to solve. This positive behaviour is consistent with its design and computational results for pure BEM formulations (cf.~\cite{antoine2021introduction, wout2021pmchwt}). Finally, an interesting observation is that OSRC preconditioning also works well when applied to the pure FEM operator in the standard FEM-BEM system~\eqref{eq:system:standard}.

\FloatBarrier
\subsubsection{ILU preconditioning}

The operator preconditioning strategies only apply to the BEM operators in the system matrix, not to the FEM operator. Remember that the top left block of the system matrix includes the summation of a FEM and BEM operator, where the boundary integral operator only acts on the boundary elements of the volumetric variational formulation (cf.~Eq.~\eqref{eq:system:stable}). As explained in Section~\ref{sec:prec:fembem}, two different preconditioning strategies can be distinguished: 1) using the ILU preconditioner for all degrees of freedom, including the boundary nodes; or, 2) using the ILU preconditioner for the inner nodes only and an operator preconditioner for the boundary nodes. 

\begin{figure}[!ht]
	\centering
	\includegraphics[width=\columnwidth]{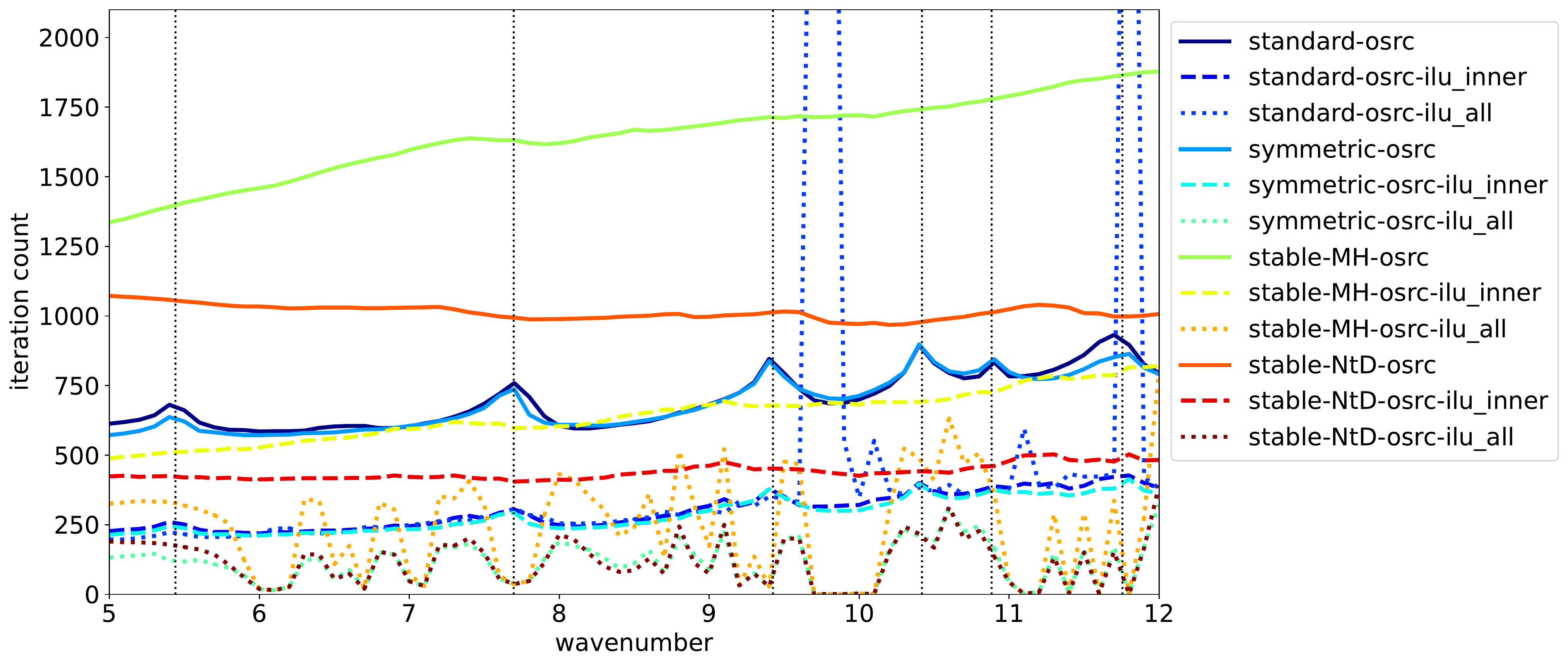}
	\caption{The influence of ILU preconditioning on the performance of the different FEM-BEM coupling techniques. The stabilisation parameters are $\eta=1$ and $\nu=0$, P1 elements are used, and OSRC preconditioning is performed.}
	\label{fig:cube:ilu}
\end{figure}

The results in Figure~\ref{fig:cube:ilu} show that both versions of FEM-BEM coupled preconditioning improve GMRES convergence significantly. Generally speaking, using all FEM nodes for the ILU preconditioner is more efficient but less robust, showing an erratic behaviour with respect to frequency. At several simulations, the ILU preconditioner invalidates an effective stabilisation of the FEM-BEM system which could be due to an imbalance of the stabilisation term in the preconditioned norm. Dedicated ILU preconditioning strategies might alleviate these issues but this study only considers readily available implementations with default settings as preconditioning. Furthermore, allowing for more fill-in elements in the ILU preconditioner or using a full LU factorisation for the FEM block improves convergence but at the expense of a higher computational complexity.

\FloatBarrier
\subsubsection{Row permutation and inverse regularisation}

The operators on the diagonal of a block system typically have a stronger influence on the GMRES convergence than the off-diagonal blocks. This observation is exploited by the preconditioning strategy of permuting block rows in the system matrix. Notice that using the alternative formulation~\eqref{eq:system:stable:alternative:reg} will have a similar effect since it changes the operator that is present on the third diagonal block.

For brevity, let us investigate the influence of row permutation and inverse regularisation for the most efficient formulations. The previous benchmarks confirm that the OSRC-NtD regulariser yields lower GMRES iteration counts than the modified Helmholtz and shifted Laplacian as stabilisation operator. Furthermore, the combination of OSRC preconditioning for the BEM operators and ILU preconditioning of the inner FEM nodes improves the GMRES convergence considerably while remaining robust.

\begin{figure}[!ht]
	\centering
	\includegraphics[width=\columnwidth]{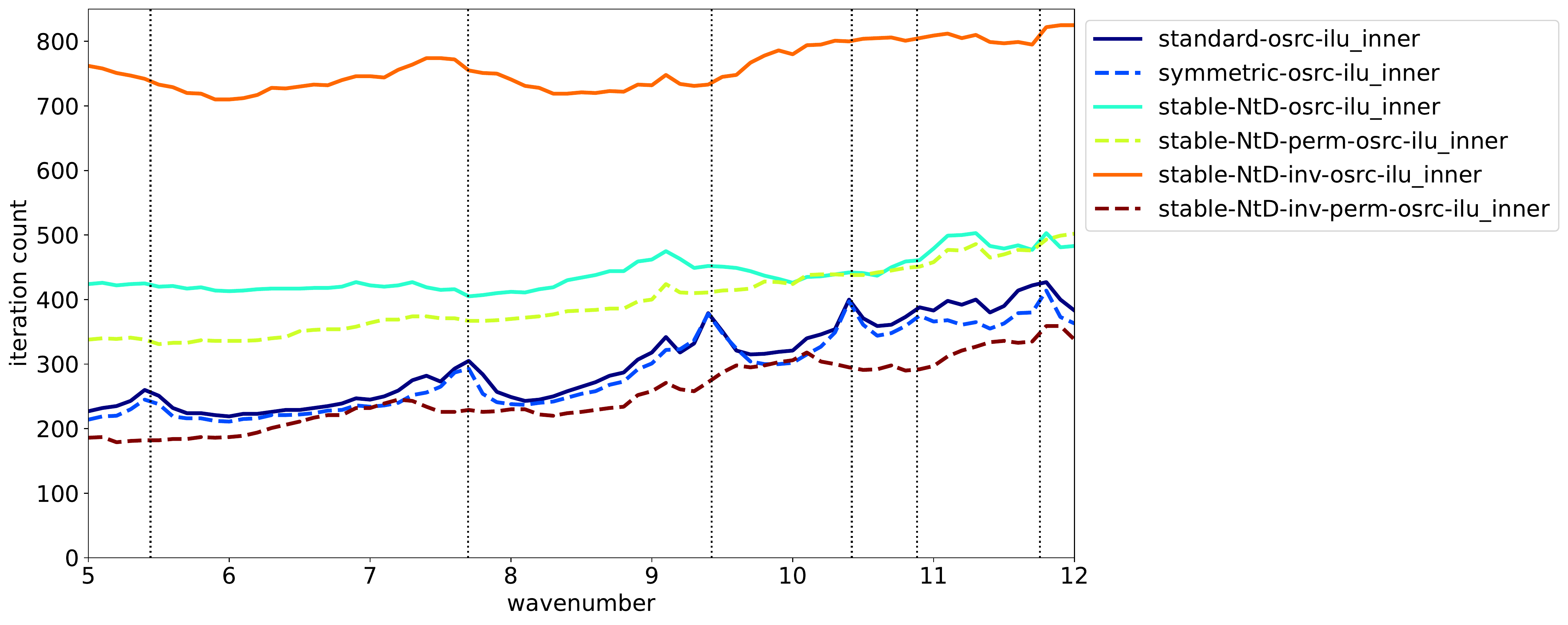}
	\caption{The influence of row permutation and inverse regularisation on the performance of the NtD-regularised FEM-BEM coupling techniques. The stabilisation parameters are $\eta=1$ and $\nu=0$, P1 elements are used, and combined OSRC-ILU preconditioning is used.}
\label{fig:cube:ntd}
\end{figure}

Figure~\ref{fig:cube:ntd} confirms that simple row permutation can be an effective preconditioning strategy depending on the operators in the system. Also, taking the inverse preconditioner changes the system in two ways: it includes the NtD operator instead of the DtN operator and the stabilisation operator is present at a different location in the block matrix. In the end, the formulation that leads to quickest GMRES convergence in this benchmark is the alternative formulation~\eqref{eq:system:stable:alternative:reg} with row permutation.

\FloatBarrier
\subsubsection{Benchmark summary}

The computational results show that the OSRC-NtD regulariser yields the stabilised FEM-BEM formulation with fastest GMRES convergence. A Galerkin discretisation with P1 elements is more efficient than a mixed P0-P1 scheme. Furthermore, OSRC preconditioning reduces the iteration count and can be combined with an ILU preconditioner for the inner nodes. Finally, minor improvements can be achieved by considering the alternative FEM-BEM formulations.

\FloatBarrier
\subsection{Large-scale simulation}

As a final benchmark to showcase the performance of the FEM-BEM methodology, a large-scale simulation will be performed. The geometry will not only consider a heterogeneous wavenumber in the interior but also a heterogeneous density and multiple domains. The heterogeneous density in the interior influences the definition of the FEM matrix as well as the density ratios in the FEM-BEM coupled system. The stabilised FEM-BEM algorithm can be extended straightforwardly to multiple domains where the cross interactions between domains are modelled through boundary integral operators, see~\ref{sec:multiple}.

\begin{figure}[!ht]
	\centering
	\begin{subfigure}[b]{\columnwidth}
		\centering
		\includegraphics[width=.8\columnwidth]{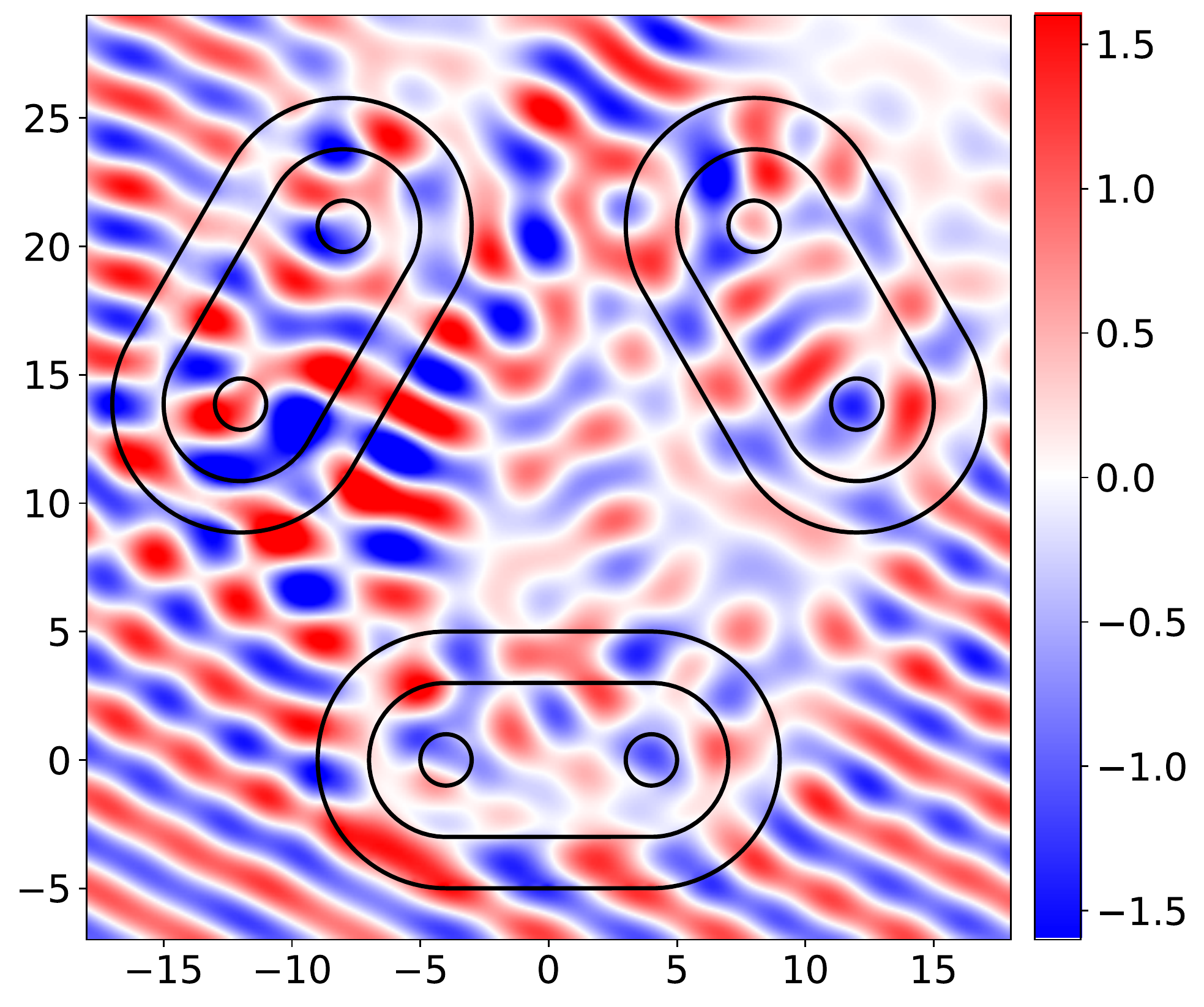}
		\caption{The real part of the acoustic field at the plane $z=0$.}
	\end{subfigure}
	\begin{subfigure}[b]{\columnwidth}
		\centering
		\vspace*{5mm}
		\includegraphics[width=.95\columnwidth]{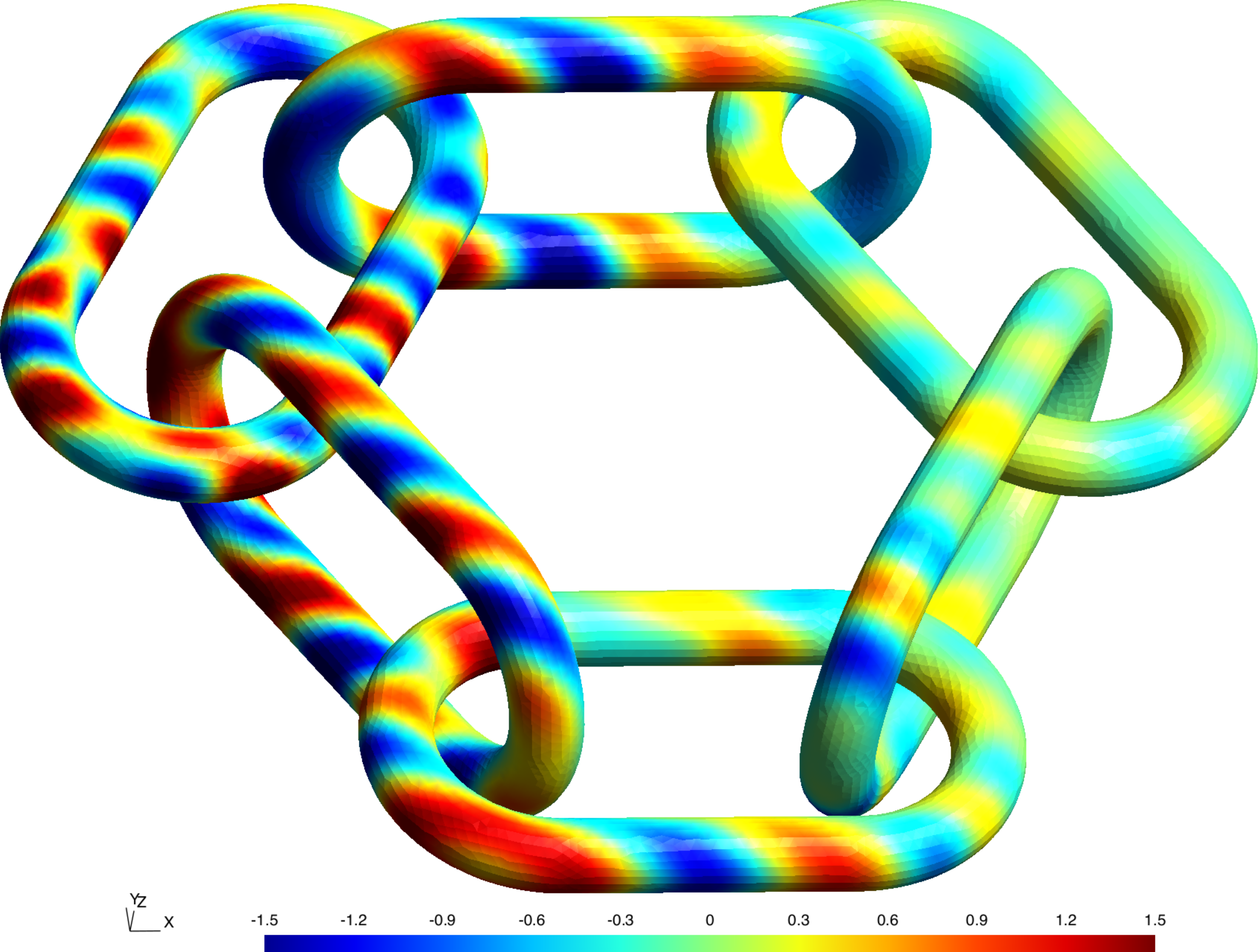}
		\caption{The real part of the acoustic field at the surface of each object.}
	\end{subfigure}
	\caption{The acoustic field calculated by the OSRC-stabilised FEM-BEM algorithm with combined OSRC-ILU preconditioning.}
	\label{fig:chain}
\end{figure}

As computational domain, let us consider a chain of six interlocked rings. Each ring has a (nondimensionalised) thickness of 2, the curves have a radius of 4 to the centre, and the straight tubes have a length of 4. Three of the rings are oriented horizontally and three vertically, see Figure~\ref{fig:chain}. The exterior wavenumber is 2 so that the exterior wavelength is $\pi$. The maximum cell size is 0.64 so that at least 5 elements fit within a wavelength. In total, there are 36\,034 nodes in the six volumetric meshes and 21\,770 nodes in the surface meshes. This yields a linear system with 57\,804 degrees of freedom. As refractivity, the functions
\begin{subequations}
\begin{align}
	n_1(\mathbf{x}) &= 0.5\frac{18-x}{36}+0.1\frac{x+18}{36}, \\
	n_2(\mathbf{x}) &= 0.5-0.4\frac{25-z^2}{25}
\end{align}
\end{subequations}
are used for the horizontally and vertically oriented rings, respectively. The interior density was chosen as
\begin{equation}
	\rhoint(\mathbf{x}) = 0.2+5\left(\frac{x+18}{36}\right)^2.
\end{equation}
Here $\mathbf{x} = (x,y,z)$ denotes the position. Notice that the wave speed and density are smooth in the interior but have a jump at the material interface. The incident field is a plane wave with direction $(1/\sqrt{5},2/\sqrt{5},0)$. The FEM-BEM algorithm uses the stabilised formulation that considers OSRC regularisation and combined OSRC-ILU preconditioning. The discretisation uses P1 functions. The same numerical parameters as specified in Section~\ref{sec:results:settings} were used, for example, a GMRES tolerance of $10^{-5}$. The GMRES solver converged in 1781 iterations and Figure~\ref{fig:chain} shows the acoustic field.

\section{Conclusions}

The standard FEM-BEM algorithms require stabilisation procedures to avoid spurious solutions at resonance frequencies of the structure. The regularisation of the surface potential guarantees robustness but deteriorates the conditioning of the linear system. This study proposed two novelties to improve the convergence of the GMRES linear solver. Firstly, the OSRC approximation of the NtD operator is used as the regulariser for the stabilisation term. Secondly, a combined preconditioning strategy of OSRC operators for the BEM blocks and ILU factorisation for the inner FEM degrees of freedom significantly reduces the number of GMRES iterations for all FEM-BEM formulations considered in this study. Computational benchmarks on canonical test cases showcase the performance of the proposed methodology. Finally, a large-scale simulation of acoustic propagation through six interlocked rings was presented. The preconditioned FEM-BEM algorithm accurately simulates the wavefield in the entire region, where the material had heterogeneous parameters for both the wave speed and density. Additional gains can be achieved in future research by considering fast algorithms for the dense matrix arithmetic and tuning the numerical parameters in the OSRC operators and ILU factorisation.

\FloatBarrier
\appendix

\section{Calderón projections}
\label{sec:calderon}

The incident, scattered and total acoustic pressure fields $\pinc$, $\psca$ and $\ptot$ all satisfy the Helmholtz equation with wavenumber $\kext$ in the homogeneous, unbounded exterior domain $\Omega^+$. Furthermore, the scattered field satisfies the Sommerfeld radiation condition, see Eq.~\eqref{eq:helmholtz}. Hence, the scattered field satisfies the representation formula~\eqref{eq:representation}, that is,
\begin{equation*}
	\psca = \SLP(\psi) - \DLP(\phi)
\end{equation*}
for the surface potentials $\phi$ and $\psi$ and potential operators $\SLP$ and $\DLP$. The traces of the potential operators satisfy the jump relations (cf.~\cite{nedelec2001acoustic, steinbach2008numerical, sauter2010boundary}):
\begin{subequations}
\begin{align}
	\traceDe(\SLP\psi) &= \SL\psi,
	& \traceDi(\SLP\psi) &= \SL\psi, \\
	\traceDe(\DLP\phi) &= \DL\phi + \tfrac12\phi,
	& \traceDi(\DLP\phi) &= \DL\phi - \tfrac12\phi, \\
	\traceNe(\SLP\psi) &= \AD\psi - \tfrac12\psi,
	& \traceNi(\SLP\psi) &= \AD\psi + \tfrac12\psi, \\
	\traceNe(\DLP\phi) &= -\HS\phi,
	& \traceNi(\DLP\phi) &= -\HS\phi,
\end{align}
\end{subequations}
where $\SL$, $\DL$, $\AD$, and $\HS$ are the boundary integral operators~\eqref{eq:boundaryoperators}. The traces of the representation formula read
\begin{equation}
	\label{eq:representation:surface}
	\begin{bmatrix} \traceDei\psca \\ \traceNei\psca \end{bmatrix} = \begin{bmatrix} -\DL \mp \tfrac12\ID & \SL \\ \HS & \AD \mp \tfrac12\ID \end{bmatrix} \begin{bmatrix} \phi \\ \psi \end{bmatrix}.
\end{equation}
Notice that
\begin{equation}
	\begin{bmatrix} \traceDi\psca \\ \traceNi\psca \end{bmatrix} - \begin{bmatrix} \traceDe\psca \\ \traceNe\psca \end{bmatrix} = \begin{bmatrix} \phi \\ \psi \end{bmatrix}.
\end{equation}
Hence, when extending the scattered field towards the interior domain, one can either use
\begin{equation}
	\label{eq:psca:1}
	\psca = \begin{cases} \ptot - \pinc, & \Omega^+; \\ -\pinc, & \Omega^-; \end{cases}
\end{equation}
leading to
\begin{equation}
	\begin{cases} \phi = -\traceDi\pinc - \traceDe\ptot + \traceDe\pinc = -\traceDe\ptot, \\ \psi = -\traceNi\pinc - \traceNe\ptot + \traceNe\pinc = -\traceDe\ptot \end{cases}
\end{equation}
or use
\begin{equation}
	\label{eq:psca:2}
	\psca = \begin{cases} \ptot - \pinc, & \Omega^+; \\ 0, & \Omega^-; \end{cases}
\end{equation}
leading to
\begin{equation}
	\begin{cases} \phi = -\traceDe\ptot + \traceDe\pinc = -\traceDe\psca, \\ \psi = -\traceNe\ptot + \traceNe\pinc = -\traceDe\psca. \end{cases}
\end{equation}
Substituting the choice of the scattered field~\eqref{eq:psca:1} or~\eqref{eq:psca:2} into the exterior boundary integral representation~\eqref{eq:representation:surface} yields either
\begin{equation}
	\begin{bmatrix} -\DL + \tfrac12\ID & \SL \\ \HS & \AD + \tfrac12\ID \end{bmatrix} \begin{bmatrix} \traceDe\ptot \\ \traceNe\ptot \end{bmatrix} = \begin{bmatrix} \traceDe\pinc \\ \traceNe\pinc \end{bmatrix}
\end{equation}
or
\begin{equation}
	\begin{bmatrix} -\DL + \tfrac12\ID & \SL \\ \HS & \AD + \tfrac12\ID \end{bmatrix} \begin{bmatrix} \traceDe\psca \\ \traceNe\psca \end{bmatrix} = \begin{bmatrix} 0 \\ 0 \end{bmatrix}.
\end{equation}
Notice that the last version can be written into the projector form
\begin{equation*}
	\begin{bmatrix} \traceDe\psca \\ \traceNe\psca \end{bmatrix} = P \begin{bmatrix} \traceDe\psca \\ \traceNe\psca \end{bmatrix}
\end{equation*}
with
\begin{equation}
	\label{eq:calderon:projector}
	P = \begin{bmatrix} \DL + \tfrac12\ID & -\SL \\ -\HS & -\AD + \tfrac12\ID \end{bmatrix} = \tfrac12\ID - A
\end{equation}
where $A$ denotes the Calderón operator.

\section{Dirichlet-to-Neumann maps}
\label{sec:dtn}

Writing out the general Calderón projection relation
\begin{subequations}
\begin{align}
	\begin{bmatrix} \traceDe\psca \\ \traceNe\psca \end{bmatrix} &= \mathbb{P} \begin{bmatrix} \traceDe\psca \\ \traceNe\psca \end{bmatrix}, \\
	\mathbb{P} &= \begin{bmatrix} \mathbb{A} & \mathbb{B} \\ \mathbb{C} & \mathbb{D} \end{bmatrix}
\end{align}
\end{subequations}
yields
\begin{align*}
	\traceDe\psca = \mathbb{A} \traceDe\psca + \mathbb{B} \traceNe\psca, \\
	\traceNe\psca = \mathbb{C} \traceDe\psca + \mathbb{D} \traceNe\psca.
\end{align*}
Hence,
\begin{align*}
	\left(\ID - \mathbb{A}\right) \traceDe\psca &= \mathbb{B} \traceNe\psca, \\
	\mathbb{C} \traceDe\psca &= \left(\ID - \mathbb{D}\right) \traceNe\psca, \\
	\mathbb{C} \traceDe\psca &= \traceNe\psca - \mathbb{D} \mathbb{B}^{-1} \left(\ID - \mathbb{A}\right) \traceDe\psca.
\end{align*}
Defining the Dirichlet-to-Neumann map $\DtNe$ implicitly as
\begin{equation}
	\DtNe \traceDe\psca = \traceNe\psca
\end{equation}
one has
\begin{subequations}
\begin{align}
	\DtNe &= \mathbb{B}^{-1} \left(\ID - \mathbb{A}\right), \\
	\DtNe &= \left(\ID - \mathbb{D}\right)^{-1} \mathbb{C}, \\
	\DtNe &= \mathbb{C} + \mathbb{D} \mathbb{B}^{-1} \left(\ID - \mathbb{A}\right).
\end{align}
\end{subequations}
Substituting the standard exterior Calderón operators into the projectors yields
\begin{subequations}
	\begin{align}
		\DtNe &= \SL^{-1} \left(\DL - \tfrac12\ID\right), \\
		\DtNe &= -\left(\tfrac12\ID + \AD\right)^{-1} \HS, \\
		\DtNe &= -\HS + \left(\tfrac12\ID - \AD\right) \SL^{-1} \left(\DL - \tfrac12\ID\right)
	\end{align}
\end{subequations}
three equivalent definitions of the DtN map.

\section{Multiple domains}
\label{sec:multiple}

Let us consider a geometry with two disjoint objects $\Omega_1$ and $\Omega_2$, and an unbounded exterior region $\Omega_0$. In the interior domain $\Omega^- = \Omega_1 \cup \Omega_2$ the FEM can be applied to the Helmholtz equation as before, with the difference that the FEM block can be decoupled into two independent blocks in the system matrix. Concerning the unbounded exterior domain, one could consider the BEM on the entire material interface $\Gamma = \partial \Omega^- = \partial \Omega_1 \cup \partial \Omega_2$ but it is more convenient and insightful to consider the two interfaces $\Gamma_2 = \partial \Omega_2$ and $\Gamma_1 = \partial \Omega_2$ separately. The stabilised FEM-BEM system~\eqref{eq:system:stable} for two scatterers reads
\begin{align*}
	\label{eq:system:stable:multiple}
	&\begin{bmatrix}
		\mathcal{F}_1 + \frac\rhoext\rhoint \HS_{11} & \frac\rhoext\rhoint \left(\AD_{11} - \tfrac12\ID_1\right) & 0 & \frac\rhoext\rhoint \HS_{12} & \frac\rhoext\rhoint \AD_{12} & 0 \\
		\tfrac12\ID_1 - \DL_{11} & \SL_{11} & \imath\eta\ID_1 & -\DL_{12} & \SL_{12} & 0 \\
		-\HS_{11} & -(\tfrac12\ID_1 + \AD_{11}) & \mathcal{S}_1 & -\HS_{12} & -\AD_{12} & 0 \\
		\frac\rhoext\rhoint \HS_{21} & \frac\rhoext\rhoint \AD_{21} & 0 & \mathcal{F}_2 + \frac\rhoext\rhoint \HS_{22} & \frac\rhoext\rhoint \left(\AD_{22} - \tfrac12\ID_2\right) & 0 \\
		-\DL_{21} & \SL_{21} & 0 & \tfrac12\ID_2 - \DL_{22} & \SL_{22} & \imath\eta\ID_2 \\
		-\HS_{21} & -\AD_{21} & 0 & -\HS_{22} & -(\tfrac12\ID_2 + \AD_{22}) & \mathcal{S}_2
	\end{bmatrix}
	\begin{bmatrix} p_1 \\ \theta_1 \\ \Sigma_1 \\ p_2 \\ \theta_2 \\ \Sigma_2 \end{bmatrix} \\
	&= \begin{bmatrix}
		\frac\rhoext\rhoint \HS_{11} \gamma_{D,1}^+ \pinc + \frac\rhoext\rhoint \gamma_{N,1}^+ \pinc \\
		(\tfrac12\ID - \DL) \traceDe\pinc \\
		-\HS \traceDe\pinc \\
		\frac\rhoext\rhoint \HS_{22} \gamma_{D,2}^+ \pinc + \frac\rhoext\rhoint \gamma_{N,2}^+ \pinc \\
		(\tfrac12\ID_2 - \DL_{22}) \traceDe\pinc \\
		-\HS_{22} \traceDe\pinc \\
	\end{bmatrix}
\end{align*}
where $\nu=0$ was taken. Here, the subindices 1 and 2 denote the operator or variable corresponding to the objects $\Omega_1$ and $\Omega_2$ or the boundaries $\Gamma_1$ and $\Gamma_2$, respectively. The double subindices indicate cross interactions, for example the operator $\SL_{mn}$ has domain $\Gamma_n$ and range $\Gamma_m$. The system can be extended to an arbitrary number of scatterers.

\bibliographystyle{unsrt}
\bibliography{refs}

\end{document}